\documentclass[reqno,11pt]{amsart}

\usepackage{amssymb,amsthm}
\usepackage{amsmath}
\usepackage{amsfonts}
\usepackage{dsfont}
\usepackage{lineno}
\usepackage{graphicx}
\usepackage[a4paper,  margin=3.2cm]{geometry}

\usepackage[active]{srcltx}
\makeatletter\@addtoreset{equation}{section}\makeatother

\newtheorem{theorem}{Theorem}[section]

\newtheorem{proposition}[theorem]{Proposition}

\numberwithin{equation}{section}

\title{On the statistical mechanics of large
populations}

\author{ Yuri  Kozitsky}
\address{Instytut Matematyki, Uniwersytet Marii Curie-Sk{\l}odowskiej, 20-031 Lublin, Poland}
\email{jkozi@hektor.umcs.lublin.pl}

\author{Krzysztof Pilorz}
\address{Instytut Matematyki, Uniwersytet Marii Curie-Sk{\l}odowskiej, 20-031 Lublin, Poland}
\email{krzysztof.pilorz@umcs.pl}

\begin{document}


\begin{abstract}

There exists a wide variety of works on the dynamics of large
populations ranging from simple heuristic modeling to those based on
advanced computer supported methods. Their interconnections,
however, remain mostly vague, which significantly limits the
effectiveness of using computer methods in this domain.  The aim of
the present publication is to propose a concept based on the
experience elaborated in the nonequilibrium statistical mechanics of
interacting physical particles. Its key aspect is to explicitly
describe micro-states of populations of interacting entities as
probability measures and then to link this description to its
macroscopic counterpart based on kinetic-like equations, suitable
for solving numerically. The pivotal notion introduced here is a
sub-Poissonian state where the large $n$ asymptotic of the
probability of finding $n$ particles in a given vessel is similar to
that for noninteracting entities, for which macro- and microscopic
descriptions are equivalent. To illustrate the concept, an
individual based model of an infinite population of interacting
entities is proposed and analyzed. For this population, its
evolution preserves sub-Poissonian states, that allows one to
describe it through the correlation functions of such states for
which a chain of evolution equations is obtained. The corresponding
kinetic equation is derived and numerically solved and analyzed.

\end{abstract}

\maketitle

\section{Introduction}

There exists a broad spectrum of publications dedicated to studying
dynamics of large populations at various levels of mathematical
sophistication, see, e.g.,
\cite{Neuhauser,BP3,Burini,otso,Cornell,OSF}. Mostly, such
populations are complex systems in the sense that their collective
behavior does not follow directly from the behaviors of the
individual parts, cf. \cite{Newman}. At the same time, a large
system becomes complex if it has an inner structure determined by
the interactions of the constituents with each other. In particular,
according to \cite{Green} ``Patterns and processes resulting from
interactions between individuals, populations, species and
communities in landscapes are the core topic of ecology. These
interactions form complex networks, which are the subject of intense
research in complexity theory, informatics and statistical
mechanics".

Before 1970s, ecological modelers employed mostly macroscopic
phenomenological (rather heuristic) models: Verhulst,
Lotka-Volterra, SIS, SEIR, etc., described by simple equations, see
\cite{Neuhauser}, which resemble phenomenological thermodynamic
models of physical substances like gases, liquids etc. In such
models -- similarly as in their physical prototypes -- the main
(often the sole) quantitative characteristic is the population
density, or a similar aggregate parameter. However, from the
statistical physics perspective one might see that this kind of
modeling provides only a coarse-grained (macroscopic) sketch of the
picture that can be seen at a microscopic level by employing
individual based models \cite{DeA}. In such models, the agents
constituting the system are explicitly taken into consideration by
assigning them individual traits, such as spatial location, type,
etc. Then interactions in the system, that form its local structure,
are described as functions of agent's traits. Therefore, similarly
as in the case of physical substances, a comprehensive description
of the dynamics of a large population should be done at at least two
scales: macro- and micro-scale. This also means that the population
itself should be considered as infinite if one wants to clearly
distinguish between the global and the local aspects of its
evolution. This seems to be true, see \cite{KozF}, even if the
population consists only of thousands -- not billions -- of members.
In statistical physics, the necessity of taking an infinite volume
(in fact, infinite particle) limit for describing phase transitions
(interaction-induced collective phenomena) was accepted after a
heated discussion followed by voting, see \cite[pages 5,6]{Simon}.
Since that time, it has become clear that the individual-based
mathematical theory of a statistically large complex system should
employ a probabilistic/statistical apparatus -- even if the
elementary evolutionary acts are deterministic.

Generally speaking, our aim in this work is to clarify the
interconnections between the macro- and the microscopic descriptions
of the evolution of large populations based on the experience
elaborated in the nonequilibrium statistical mechanics of infinite
systems of interacting physical particles and some recent
developments in this direction. This seems to be even more important
in view of the wide variety of powerful computational methods which
can efficiently be employed here. In statistical physics, the
conceptual fundamentals were laid by N. N. Bogoliubov \cite{Bogol}
(see also \cite{DSS} for a contemporary interpretation of
Bogoliubov's approach), who proposed to describe micro-states of
large systems of physical point particles as probability measures on
appropriate phase spaces, characterized by their correlation
functions. In recent years, this approach is being developed also in
the theory of large populations, see
\cite{BP3,otso,Cornell,OSF,KozF,Dima,KoKoz,Koz,OK} and the works
quoted therein. The key feature of these and similar works consists
in describing the evolution of the population states as that of the
corresponding correlation functions. Among them there are the first
and the second order functions,  describing the population density
and the binary spatial correlations, respectively. Here, similarly
as in Bogoliubov's theory, the evolution of the whole set of
correlation functions is usually governed by an infinite chain of
coupled equations, which turns their simultaneous solving into a
hard task. Approximate solutions of this chain are usually obtained
by its decoupling with the help of a certain \emph{ansatz}, which
yields closed evolution equations for the densities alone, see e.g.,
\cite{Dima}, or jointly with the binary spatial correlations, see
\cite{Cornell,OSF,OK}. An important advantage of this is that here
one can apply appropriate numerical methods, see
\cite{OK,Omel1,KOP}, to get a deeper insight into the properties of
their solutions. In these and similar works, however, it is tacitly
assumed (believed) that the description obtained in this way
reliably approximates the evolution of states of the considered
system. One of the key aims of this article is to demonstrate that
this is not always the case and that the use of low order
correlation functions only ought to be preceded by a detailed
analyzes of the evolution of states of the underlying individual
based model. To distinguish between the cases where one can or
cannot restrict the theory to studying only the low order
correlation functions, we introduce the notion of a
\emph{sub-Poissonian} state. In some sense, such states are nearly
Poissonian, and thus can sufficiently be characterized by their
correlation functions of low orders only. At the same time, the set
of such states is big enough as it includes states of thermal
equilibrium of interacting physical particles that undergo a phase
transition, see \cite{Ruelle,Rue}. Recall that a Poisson state,
corresponding to the lack of interactions, is fully characterized by
its first-order correlation function. Thus, in order to base the
description on the low order correlation functions one should
demonstrate that the evolution of states preserves the set of
sub-Poissonian states, which is one of the principal challenges of
the theory. If this is the case, then one can pass to the evolution
of correlation functions and approximate it either by decoupling the
second correlation function into the product of two densities
(mean-field approximation), or go beyond by decoupling, e.g., third
order function.

As mentioned above, we intend to use the experience elaborated in
nonequilibrium statistical physics to reconsider the steps leading
from the macro- to the microscopic dynamics of large populations. To
this end, we introduce and study a generalization of the known
Bolker-Pacala model \cite{BP1,Dima} in which the population members
at random arrive in and depart from a continuous habitat. The
arrival and departure rates depend also on the existing population.
We begin by a phenomenologically deduced evolution equation, see
(\ref{1}), subsequently transformed into the kinetic equation
(\ref{5}) in the spirit of the celebrated mean-field approach
\cite{Simon}. The latter step is preceded by the introduction of the
notions of \emph{occupation probabilities} and \emph{expected number
of $k$-particle clusters}, essential for the further analyzes. Then
the passage to the corresponding micro-model and its theory is
performed in the way that stresses the connections with Bogoliubov's
theory. To do this, we make precise the very notion of a state of a
possibly infinite population made in Sect. \ref{S2}, where we pay
attention to another essential notion: a \emph{sub-Poissonian}
state, for which the large $n$ asymptotic of the occupation
probabilities should be at most such as in a Poisson state. In
subsect. \ref{SS23}, we provide an example of the state in which the
occupation probabilities have a \emph{heavy tail}, and thus the high
order correlations essentially determine its properties. In Sect.
\ref{S3}, we introduce our microscopic model and demonstrate that
its dynamics preserves the set of sub-Poisonian states, and thus the
evolution of its states can reasonably be approximated by
kinetic-like equations. In subsect. \ref{SS34}, we present the
results of the numerical study of such equations. Finally, in Sect.
\ref{SS35} we discuss its main aspects and make some comparison,
mostly with the approach of \cite{Cornell}.

\section{Posing}


 By $\mathbf{N}$ we denote the set of positive integers
$1,2,\dots$, whereas $\mathbf{R}^d$, $d\geq 1$, will stand for the
usual space of vectors $x= (x^1, \dots x^d)$ with real components
$x^j$. It will serve us as the habitat of the considered systems. A
`vessel', $\Lambda$, is a bounded closed subset of the habitat. It
will always be involved when discussing local properties. For an
appropriate function $f$, we write
\begin{eqnarray}
  \label{0}
  \int f(x) d x & = & \int_{-\infty}^{+\infty} \cdots
  \int_{-\infty}^{+\infty} f(x^1,\dots , x^d) d x^1 \cdots d
  x^d,\\[.2cm] \nonumber
\int_\Lambda f(x) d x & = & \int f(x) \mathbf{1}_\Lambda (x) d x,
\end{eqnarray}
where
\begin{equation}
  \label{18a}
\mathbf{1}_\Lambda (x) = \left\{\begin{array}{ll} 1, \qquad &{\rm
if}  \ \ x\in  \Lambda; \\[.2cm] 0, \qquad &{\rm
otherwise}.
\end{array} \right.
\end{equation}
Sometimes, by $\int$ we denote also other integrals, the meaning of
which becomes clear from the context.

\subsection{From macro to micro}


Assume that we are given a continuous medium filling in the whole
space $\mathbf{R}^d$. Its only quantitative characteristic is given
by the density $\varrho(x)$, which means that the evolution of this
substance is described as a map $\varrho_0 (x) \to \varrho_t(x)$,
where $t\geq 0$ is time and $\varrho_0 (x)$ is an initial condition.
Assume further that it is obtained by solving the following
differential equation
\begin{equation}
  \label{1} \frac{\partial}{\partial t}\varrho_t(x) = b^{+}(x) -
  b^{-}(x) \varrho_t(x), \qquad \varrho_t(x)|_{t=0}=\varrho_0(x),
\end{equation}
where $b^{\pm}(x)\geq 0$ are certain functions, characterizing this
evolution. Their dependence on $x$ reflects possible spatial
heterogeneity of the environment. The key feature of (\ref{1}) is
that it is local in space, i.e., densities at distinct points evolve
independently of each other. To imagine the object corresponding to
this equation, one may think of a layer of a viscose liquid the
density of which changes in time due to condensation and evaporation
described by $b^{+}(x)$ and $b^{-}(x)$, respectively. The solution
of (\ref{1}) is obtained in the form
\begin{equation}
  \label{2}
  \varrho_t (x) = \varrho_0(x) e^{-b^{-}(x)t} + \left(1-e^{-b^{-}(x)t}\right)
  b^{+}(x)/b^{-}(x), \qquad b^{-}(x) >0,
\end{equation}
and $\varrho_t(x) = \varrho_0(x) + b^{+}(x) t$ for $b^{-}(x)=0$.
Note that $\varrho_t(x)$ is increasing ad infinitum in the latter
case, whereas for $b^{-}(x) >0$, $\varrho_t (x) \to
b^{+}(x)/b^{-}(x)$ as $t\to +\infty$. If the environment is
homogeneous, $b^{+}(x)/b^{-}(x)$ is independent of $x$; hence, the
spatial distribution of the substance gets asymptotically
homogeneous in this case.

Let us assume now that, instead of a continuous medium, we deal with
point entities -- particles -- arriving in and departing from the
habitat $\mathbf{R}^d$, and that $\varrho_t(x)$ is their density at
time $t$. In this case, however, $\langle
\varrho_t\rangle_\Lambda=\int_\Lambda \varrho_t(x) d x$ should be
interpreted as the expected number of particles contained in a given
vessel $\Lambda$ at time $t$, corresponding to a certain probability
law, which should now explicitly appear in the theory  as the system
state. In such a state, the system gets infinite with probability
one if $\int \varrho_t(x) d x =+\infty$. By making this change of
interpretation one passes to the microscopic scale, and therefore
starts using statistical/probabilistic tools, which is inevitable in
view of the system size. The explicit definition of the mentioned
probability law -- system's state -- could be performed upon
introducing appropriate mathematical framework. What might be
interesting for a less mathematically inclined ecological modeler is
the occupation probability $p_t(n,\Lambda)$ corresponding to this
state, i.e., the probability of finding $n$ particles contained in
$\Lambda$  at time $t$. To stay in the approach based on the use of
densities only, one should assume that there is no interaction
between the particles -- in agreement with the spacial locality of
(\ref{1}). In probability theory, the lack of interactions
corresponds to the lack of dependencies. In this case, the
aforementioned state is given by the Poisson distribution according
to which the occupation probabilities are
\begin{equation}
  \label{3}
  p_t(n,\Lambda) = \frac{\langle \varrho_t\rangle^n_\Lambda}{n!}
  \exp\left( - \langle \varrho_t\rangle_\Lambda\right).
\end{equation}
A similar characteristic of a system state is the expected number of
$k$-element clusters contained in $\Lambda$. In the Poissonian state
corresponding to (\ref{3}), it is
\begin{equation}
  \label{4}
  c_t(k, \Lambda) = [c_t(1, \Lambda)]^k= \langle \varrho_t\rangle^k_\Lambda,
  \qquad k\in \mathbf{N},
\end{equation}
which reflects the independence of particles in taking their
positions in such states. Thus, without interactions the microscopic
evolution $p_0(n,\Lambda) \to p_t(n,\Lambda)$, $c_0(k,\Lambda) \to
c_t(k,\Lambda)$ is obtained from the macroscopic one $\varrho_0(x)
\to \varrho_t(x)$, and vice versa. One may say that such a system is
not complex. If the evolution $\varrho_0(x) \to \varrho_t(x)$ is
governed by (\ref{1}), the time dependence of the microscopic
quantities (\ref{3}) and (\ref{4}) can be obtained explicitly, see
(\ref{2}). Note that, in the `corpuscular' interpretation, the first
term in (\ref{1}) describes the attraction of the newcomers by the
environment, whereas the second one corresponds to the disappearance
(emigration) of the existing population members, also caused by the
environment.

Let us figure out now how possible interactions between the
particles can be introduced to the theory. In statistical physics,
the simplest way of doing this without passing to the micro-scale is
known since the time of Hartree and van der Waals. Nowadays,  it is
called a \emph{naive mean-field approach}, see \cite[page
216]{Simon}. In its framework, interactions between the particles
appear as self-consistent (i.e., state-dependent) additive
`corrections' to their interactions with the environment. For the
system discussed above, such corrections might lead to the following
modification of (\ref{1})
\begin{eqnarray}
  \label{5}
\frac{\partial}{\partial t}\varrho_t(x) & = & b^{+}(x) + \int
a^{+}(x-y) \varrho_t (y) dy\\[.2cm] & - & \left[  b^{-}(x) + \int a^{-}(x-y)
\varrho_t (y) dy\right]\varrho_t(x), \nonumber
\end{eqnarray}
where $a^{\pm}(x)\geq 0$ are suitable functions. The new terms here
describe attraction and repulsion caused by the existing population.
Their presence makes (\ref{5}) nonlinear and nonlocal, and thus
essentially complicates its mathematical study.  Usually, for such
equations one proves the existence of solutions in certain classes
of functions, and then -- by means of specially elaborated
algorithms, see, e.g., \cite{Omel1} -- constructs the solutions
numerically \cite{KOP}. This study provides a more detailed and
comprehensive description of the system dynamics in the mean-field
approximation. A more advanced way of taking interactions into
account in the mean-field approximation relies on developing a
microscopic description of a finite system with Curie-Weiss type
interactions and then passing to the infinite system limit, see,
e.g., \cite{KK}.

\subsection{Further development}

The next step from  (\ref{1}), (\ref{5}) towards the microscopic
theory  can be based on the following arguments. Assume that the
probability law -- the system state -- is such that the expected
number of $k$-particle clusters $c_t(k,\Lambda)$ contained in the
vessel $\Lambda$ at time $t$ is obtained in the form
\begin{equation}
  \label{6}
  c_t (k,\Lambda) = \int_\Lambda \cdots \int_\Lambda \varrho_t^{(k)}(x_1 , \dots ,
  x_k) d x_1 \cdots d x_k, \qquad k\in \mathbf{N},
\end{equation}
where the functions $\varrho_t^{(k)}(x_1 , \dots ,
  x_k)$, $k\in \mathbf{N}$, characterize this probability law. In fact, these $\varrho_t^{(k)}$ are
  correlation functions, see \cite{Lenard,Ruelle} and below. The first member
  of the family $\{\varrho_t^{(k)}\}_{k\in \mathbf{N}}$ is the density function used above, i.e., $\varrho_t^{(1)}(x) =
  \varrho(x)$. Then the `true' version of (\ref{5}) might be
\begin{eqnarray}
  \label{7}
\frac{\partial}{\partial t}\varrho^{(1)}_t(x) & = & b^{+}(x) + \int
a^{+}(x-y) \varrho^{(1)}_t (y) dy\\[.2cm] & - &   b^{-}(x) \varrho^{(1)}_t (x) - \int a^{-}(x-y)
\varrho^{(2)}_t (x,y) dy , \nonumber
\end{eqnarray}
and the passage from (\ref{7}) to (\ref{5}) is performed by the so
called \emph{moment closure}, cf. \cite{Mu,BP1}, i.e., by the
`ansatz' (closing approximation) $\varrho^{(2)}_t (x,y) \simeq
\varrho_t (x)\varrho_t (y)$. At this point, however, it is just a
guess based on the form of the second line in (\ref{5}). Looking at
(\ref{7}) one may surmise that it is just the first member of an
infinite chain of coupled equations. The second one would have the
time derivative of $\varrho^{(2)}_t (x,y)$ on the left-hand side,
with the right-hand side containing $\varrho^{(k)}_t$ with
$k=1,2,3$, cf. \cite{OK,Plank}. By solving this chain  one could
calculate all $\varrho_t^{(k)}$ and thereby get $c_t(k,\Lambda)$,
$p_t(k,\Lambda)$, i.e., obtain information on the microscopic
evolution of the system. Note, however, that there might be
difficult to prove that the obtained solutions $\varrho_t^{(k)}$,
$k\geq 1$ are indeed correlation functions, suitable to give the
expectations as in (\ref{6}), see \cite{Lenard}. In less
mathematically rigorous works, this kind of the moment problem \cite
{ST} is typically ignored.

At this stage, it becomes clear that the origin of the construction
yielding the equations involving $\varrho^{(k)}_t$ should be based
on a microscopic model, for which one would derive the corresponding
chain of equations for $\varrho_t^{(k)}$ and then study its
solvability, properties of the solutions, and so on. For systems of
physical particles, the first step in this direction was made by N.
N. Bogoliubov around 75 years ago in his seminal monograph
\cite{Bogol}, see also \cite{DSS} for a more contemporary version of
this theory. It was suggested there to describe the states of a
large finite system of interacting particles by means of probability
density functions $f^{(n)}_{t,\Lambda}( x_1, \dots x_n)$, where $n$
is the number of particles in the system contained in $\Lambda$ and
$x_i = (q_i,p_i)$, $q_i$ and $p_i$ being the position and the
momentum of $i$-th particle, respectively. Noteworthy, these
probability density functions exist only for a finite system
dwelling in a compact $\Lambda\subset \mathbf{R}^d$. Then the
functions $\varrho_t^{(k)}$ are obtained as follows. Their versions
in the vessel $\Lambda$ are
\begin{gather}
  \label{7a}
\varrho_{t,\Lambda}^{(k)}(x_1, \dots , x_k) = \sum_{n=0}^{+\infty}
\frac{z^n}{n!} \int_{\Lambda^*} \cdots \int_{\Lambda^*}
f^{(k+n)}_{t,\Lambda} (x_1, \dots x_k, x_{k+1} , \dots ,x_{k+n} )d
x_{k+1}\cdots dx_{k+n},
\end{gather}
where $z>0$ is activity and $\Lambda^* = \Lambda \times
\mathbf{R}^d$ is the phase space of a single particle for which $q$
varies in $\Lambda$ and $p$ can take any value in $\mathbf{R}^d$.
The system evolution is governed by the equation that can
symbolically be written in the form
\begin{equation}
  \label{7b}
  \frac{\partial f_{t,\Lambda}}{\partial t} = \{H_\Lambda,f_{t,\Lambda}\},
\end{equation}
where $H_\Lambda=\{H^{(n)}_\Lambda\}_{n\geq 2}$ is the collection of
the Hamiltonians that describe a given system -- including
interactions -- and $\{\cdot, \cdot\}$ are \emph{Poisson brackets},
see see \cite[Chapt. 2]{Bogol} for more detail. In (\ref{7b}) --
derived from Newton's laws --
$f_{t,\Lambda}=\{f^{(n)}_{t,\Lambda}\}_{n\geq 2}$ is the collection
of all probability distribution densities, which means that the full
version of (\ref{7b}) is an infinite chain of equations, from which
the equations for $\varrho_{t,\Lambda}^{(k)}$ are derived. Its
infinite volume limit $\Lambda \to \mathds{R}^d$ -- known as
Bogoliubov's chain -- is obtained by a special procedure, in which
the proof of the existence of this limit is a hard problem. This
difficulty stems from the fact that the limiting probability laws
may have correlation functions only in a generalized sense (i.e., as
distributions), for which equations like (\ref{7}) do not make
sense.

Instead of solving (\ref{7b}) one may consider the corresponding
equation for test functions (observables)  -- appropriate functions
defined in the phase space. They are introduced by the condition
\begin{gather*}
\int_{\Lambda^*} \cdots \int_{\Lambda^*}  F^{(n)}_{t,\Lambda} (x_1 ,
\dots , x_n) f^{(n)}_{0,\Lambda} (x_1 , \dots , x_n) d x_1 \cdots d
x_n \\[.2cm] \nonumber = \int_{\Lambda^*} \cdots \int_{\Lambda^*}  F^{(n)}_{0,\Lambda} (x_1 ,
\dots , x_n) f^{(n)}_{t,\Lambda} (x_1 , \dots , x_n) d x_1 \cdots d
x_n.
\end{gather*}
Then their evolution is described by the following dual version of
(\ref{7b})
\begin{equation}
  \label{7d}
  \frac{\partial F_{t,\Lambda}}{\partial t} = \{F_{t,\Lambda}, H_\Lambda\},
\end{equation}
where $F_{t,\Lambda}$ is the collection of all
$F_{t,\Lambda}^{(n)}$. A more detailed form of (\ref{7d}) is
\begin{eqnarray}
  \label{7dd}
   \frac{\partial}{\partial t}F^{(n)}_{t, \Lambda} (x_1 , \dots ,
  x_n) & = & \sum_{i=1}^n \bigg{(}\frac{\partial}{\partial q_i}F^{(n)}_{t, \Lambda} (x_1 , \dots ,
  x_n) \frac{\partial }{\partial p_i} H^{(n)}_\Lambda (x_1 , \dots ,
  x_n) \\[.2cm] & &  - \frac{\partial}{\partial p_i}F^{(n)}_{t, \Lambda} (x_1 , \dots ,
  x_n) \frac{\partial }{\partial q_i} H^{(n)}_\Lambda (x_1 , \dots ,
  x_n) \bigg{)}. \nonumber
\end{eqnarray}
States of thermal equilibrium (Gibbs states) of systems of physical
particles can also be described by probability densities, and thus
by the functions $\varrho_{\Lambda}^{(k)}$ obtained from them
according to (\ref{7a}), with the subsequent passage to the
`infinite-volume' limit $\Lambda \to \mathbf{R}^d$. For systems of
particles interacting via `super-stable potentials' \cite{Ruelle},
the infinite-volume limits of $\varrho_{\Lambda}^{(k)}$, $k\in
\mathbf{N}$ exist. They satisfy \emph{Ruelle's bound}
\begin{equation}
  \label{8}
  0 \leq \varrho^{(k)}(x_1, \dots , x_k) \leq \varkappa^k,
\end{equation}
holding for some $\varkappa>0$ and all $k$. An immediate consequence
of the latter is that, for the Gibbs states corresponding to
super-stable interaction, the expected number of $k$-clusters
 satisfies
\begin{equation}
  \label{9}
c(k,\Lambda) \leq (\varkappa |\Lambda|)^k,
\end{equation}
that resembles the case of non-interacting particles, cf. (\ref{4}).
Here $|\Lambda|$ is the volume of the vessel $\Lambda$.

For interacting particles, the $n$-particle Hamiltonians
$H_\Lambda^{(n)}$ that participate in the chains in (\ref{7b}),
(\ref{7d}) contain interaction terms. In the case of Curie-Weiss
binary interactions, they are of the form $J/n$, i.e., each two
particles interact with the same strength proportional to $1/n$. For
such interactions, the decoupling $\varrho^{(2)}_t(x,y) \simeq
\varrho^{(1)}_t(x)\varrho^{(1)}_t(y)$ is obtained rigorously in the
limit $n\to +\infty$, $\Lambda \to \mathbf{R}^d$. However, the
interaction parameter $J$ persists in the corresponding equation for
$\varrho^{(1)}_t(x)$. This procedure can be considered as a rigorous
version of the mean-field approximation mentioned above. In view of
(\ref{4}), this approximation does not change the sub-Poissonian
bound (\ref{9}) for the clustering parameter $c(k,\Lambda)$, which
points to its acceptability for states obeying this bound.

\subsection{The program}

In the remaining part of this article, we introduce and study the
microscopic model of an infinite population of point particles that
at random arrive in and depart from the habitat $\mathbf{R}^d$ in
the way corresponding to the heuristic kinetic equation (\ref{5}).
Our aim here is to follow the way similar to that of the statistical
mechanics of physical particles based on the evolution equations in
(\ref{7a}), (\ref{7b}) and (\ref{7d}). Here, however, we are going
to modify this approach as follows. Since the dynamics of our model
proves simpler than in the case of physical particles -- mostly due
to the lack of differential operators, see the right-hand side of
(\ref{7dd}) -- we can avoid problems with taking infinite volume
limits $\Lambda \to \mathbf{R}^d$ of the correlation functions
$\varrho^{(n)}_{t,\Lambda}$. This will allow us to deal with
probability laws as systems's states, `represented' by their
correlation functions that obey Ruelle's bound (\ref{8}). To this
end, we introduce the class of such laws, which we call
sub-Poissonian, and then show that the evolution in question
preserves this set. That is, the state at time $t>0$ is
sub-Poissonian if it is sub-Poissonian at $t=0$. This will allow us
to pass to the correlation functions of such states, the evolution
equation for which is derived from that governing the evolution of
states. Then we may approximate the former equations by those
containing the correlation functions of low orders only, admissible
for studying by various numerical methods.

In more detail, our plan consists of the following items.
\begin{enumerate}
\item Presenting a general framework in which one operates with
probability laws describing states of the system. Introducing
sub-Poissonian states and finding expressions for the corresponding
occupation probabilities. This is done in subsects. \ref{SS21} and
\ref{SS22}.
\item Presenting an example of the state for which the occupation
probabilities have a heavy tail, which essentially differs this
state from sub-Poissonian states. This is done in subsect.
\ref{SS23} to elucidate the importance of the latter states.
\item Introducing the model and deriving the evolution equations for
its states, similarly as in (\ref{7b}) and (\ref{7d}) for physical
systems. This is done in subsect. \ref{SS31}.
\item Checking whether this evolution preserves sub-Poissonian
character of the states, cf. (\ref{9}), which might allow one to
employ correlation functions. This is done in subsect. \ref{SS32}.
\item Deriving the corresponding chain of evolution equations for these functions, the first member of which is given in (\ref{7}).
Deriving (\ref{5}) from this chain. This is done in subsect.
\ref{SS33}.
\item Presenting the results of the numerical solving of (\ref{5}) based on the algorithm elaborated in
\cite{Omel1}. This is done in subsect. \ref{SS34}.
\end{enumerate}

\section{States of Infinite Populations}
\label{S2}

In this  section, we briefly present mathematical tools for
describing infinite populations.
\subsection{The states}

\label{SS21}

The population which we are going to study consists of point
entities -- particles -- distributed over the habitat
$\mathbf{R}^d$. Their evolution amounts to appearance (immigration)
and disappearance (emigration or death), both affected by the
existing population and the environment. From the very beginning, we
allow the population to be infinite. At the same time, we will
assume that it is  locally finite, i.e., each vessel $\Lambda$ can
contain only finitely many of its particles. Contrary to the case of
physical particles characterized by their positions and momenta $x_i
= (q_i,p_i)$, here each particle is fully characterized by its
position $x_i\in \mathbf{R}^d$ as we do not assume any spatial
motion in our model. This advantage proves decisive for constructing
solutions of the corresponding evolution equations. Therefore, pure
states of the system are infinite collections of the positions
$x_i\in \mathbf{R}^d$, which we will denote $\mathbf{x}$. By writing
$x\in \mathbf{x}$ we mean that this $x$ is a member of the
collection, i.e., $x=x_i$ for some $i$. Keeping this in mind, by
writing $\sum_{x\in \mathbf{x}} f(x)$ we will mean $\sum_{i}
f(x_i)$, where in contrast to (\ref{7dd}) the sum is infinite.
However, this latter fact makes no harm if the function $f$ has
bounded support, i.e., vanishes outside a certain vessel $\Lambda$.
Note that we allow more than one particle have the same location,
i.e., some of $x_i$'s in a given $\mathbf{x}$ may coincide. For
$x\in \mathbf{x}$,  by writing $\mathbf{x}\setminus x$ we mean the
collection obtained from $\mathbf{x}$ by removing $x$. Similarly, by
writing $\mathbf{x} \cup x$ we mean the collection obtained from
$\mathbf{x}$ by adding the mentioned $x$. The mathematical theory of
spaces consisting of such collections is well-developed, see, e.g.,
\cite{DV}.

Let $F(\mathbf{x})$ be a real-valued function, which appears in
equations like (\ref{7d}). As above, we call such functions
observable or test functions. Following Bogoliubov \cite{Bogol}, we
will consider the next class of such functions. Let $\Theta$ denote
the set of all continuous functions $\theta(x)$ such that
$\theta(x)=0$ whenever $x$ is taken outside a $\theta$-specific
vessel $\Lambda$. For $\theta \in \Theta$, we then set
\begin{equation}
  \label{13}
  F^\theta (\mathbf{x}) = \prod_{x\in \mathbf{x}} (1 + \theta (x)) = \exp\left( \sum_{x\in \mathbf{x}} \ln (1+\theta(x)) \right).
\end{equation}
If $\theta(x) \in (-1,0]$, then $0< F^\theta(\mathbf{x})\leq 1$ for
all $\mathbf{x}$. We extend (\ref{13}) also to the case where the
collection $\mathbf{x}$ is empty. This collection is denoted
$\mathbf{0}$, for which
\begin{equation}
  \label{j}
F^\theta (\mathbf{0}) = \prod_{x\in \mathbf{0}} (1 + \theta (x)) =
1.
\end{equation}
As mentioned above, we consider possibly infinite populations. This,
in particular, means that one cannot expect that the probability
laws which describe the population states have densities, which
excludes the use of formulas similar to (\ref{7a}), (\ref{7b}).
Nevertheless, such states, denoted $\mu$, are real-valued linear
maps $F\mapsto \mu(F)$ such that: (a) $\mu(F) \geq 0$ whenever
$F(\mathbf{x}) \geq 0$ for all $\mathbf{x}$; (b) $\mu(F) = C$
whenever $F(\mathbf{x}) \equiv C$. In terms of Lebesgue's integrals,
they are  $\mu(F) = \int F d \mu$, where the integrals are taken
over the set of all $\mathbf{x}$ and $d\mu$ is a probability measure
defined on this set. In the case of finite populations contained in
a given vessel $\Lambda$ and states characterized by densities, this
amounts to the following
\[
\mu(F) = F^{(0)} f_\mu^{(0)}+ \sum_{n=1}^\infty\frac{1}{n!}
\int_\Lambda \cdots \int_\Lambda F^{(n)} (x_1, \dots , x_n)
f^{(n)}_\Lambda (x_1, \dots , x_n) d x_1 \cdots d x_n,
\]
which is pretty similar to the grand canonical average, cf.
\cite{DSS}. Here $F^{(0)}=F(\mathbf{0})$ is just a real number
whereas $f_\mu^{(0)}$ is the probability that the system consists of
zero particle, i.e., is empty. Bogoliubov's brilliant idea consists
in using the following representation, cf. (\ref{j}),
\begin{equation}
  \label{14}
  \mu(F^\theta) = 1+ \sum_{n=1}^\infty\frac{1}{n!} \int \cdots \int
  \varrho_\mu^{(n)} (x_1 , \dots , x_n) \theta (x_1) \cdots \theta(x_n) d
  x_1 \cdots d x_n,
\end{equation}
that holds for all $\theta\in \Theta$ and those states $\mu$ that
have correlations functions $\varrho_\mu^{(n)}$ of all orders
$n=0,1, \dots$, such that the sum and the integrals in (\ref{14})
converge. Fortunately, the collection of such states is big enough
to develop a comprehensive theory of this kind, in which the map
$\theta \mapsto \mu(F^\theta)$ is called Bogoliubov's functional.
More on this subject can be found in
\cite{DSS,Lenard,Ruelle,Dima,KoKoz,Koz}.

\subsection{Sub-Poissonian states}
\label{SS22}

It can be shown that correlation functions of all orders are
symmetric with respect to the permutations of its arguments. That
is, for each $n\geq 2$,
\begin{equation}
  \label{sigma}
\varrho^{(n)}_\mu (x_1 , \dots , x_n) = \varrho^{(n)}_\mu
(x_{\sigma(1)} , \dots , x_{\sigma(n)}),
\end{equation}
holding for all $\sigma\in S_n$, where the latter is the symmetric
group. In particular, this means $\varrho^{(2)}_\mu (x_1,x_2) =
\varrho^{(2)}_\mu (x_2,x_1)$.

Among the states characterized by correlation functions are Poisson
states $\pi_\varrho$, for which
\begin{equation}
  \label{15}
\varrho_{\pi_\varrho}^{(n)} (x_1 , \dots , x_n) = \varrho(x_1)
\cdots \varrho(x_n), \qquad n\geq 1,
\end{equation}
where $\varrho(x)\geq 0$ is density that completely characterizes
the state $\pi_\varrho$. In this case, the occupation probabilities
$p(n,\Lambda)$ and the expectations $c(n,\Lambda)$ are as in
(\ref{3}) and (\ref{4}), respectively. By (\ref{15}) and (\ref{14})
we then get
\begin{equation}
  \label{16}
  \pi_\varrho(F^\theta) = \exp\left(\int \theta (x) \varrho (x) d x
  \right).
\end{equation}
The homogeneous Poisson state $\pi_\varkappa$ corresponds to
$\varrho(x)\equiv \varkappa>0$. In this state, the particles are
independently distributed over the space with constant density.
Keeping these facts in mind we introduce the following notion. A
state $\mu$ is called \emph{sub-Poissonian} if its Bogoliubov
functional $\mu(F^\theta)$ can be written as in (\ref{14}) and its
correlation functions $\varrho_\mu^{(k)}$ satisfy the Ruelle bound
(\ref{8}) with some $\varkappa>0$. The least such $\varkappa$ will
be called the \emph{type} of this $\mu$. For sub-Poissonian states,
the expected number of $n$-clusters in $\Lambda$ satisfies
(\ref{9}). Hence, dense clusters are not more expected in such
states than they are in the Poisson state $\pi_\varkappa$ with the
corresponding $\varkappa$. To figure out the large $n$ asymptotic of
the occupation probabilities, we proceed as follows. First we take
$\theta (x) = \zeta \mathbf{1}_\Lambda (x)$, see (\ref{18a}), with
real $\zeta$, and set
\begin{equation}
  \label{18}
\phi_{\mu,\Lambda} (\zeta) = \mu (F^{\zeta \mathbf{1}_\Lambda}).
\end{equation}
By (\ref{13}) it follows that
\[
F^{\zeta \mathbf{1}_\Lambda} (\mathbf{x}) = \prod_{x\in \mathbf{x}}
(1+ \zeta\mathbf{1}_\Lambda(x)) = (1+\zeta)^{N_\Lambda
(\mathbf{x})},
\]
where $N_\Lambda (\mathbf{x})$ is the number of particles in the
collection $\mathbf{x}$ contained in $\Lambda$. Then
\begin{gather*}
\phi_{\mu, \Lambda} (\zeta)  = \mu\left((1+\zeta)^{N_\Lambda}\right)
= \sum_{n=0}^\infty (1+\zeta)^n p_\mu(n,\Lambda), \nonumber
\end{gather*}
which yields the following formula for the occupation probabilities
\begin{equation}
  \label{18b}
p_\mu(n,\Lambda) = \phi^{(n)}_{ \mu, \Lambda} (-1)/n!,
\end{equation}
where $\phi^{(n)}_{ \mu, \Lambda} (-1)$ stands for the usual
derivatives of $\phi$ with respect to $\zeta$ at $\zeta = -1$. Note
that by (\ref{16}) it follows that $\phi_{ \pi_\varrho,
\Lambda}(\zeta) = \exp\left(\zeta \langle
\varrho\rangle_\Lambda\right)$, which agrees with the latter and
(\ref{3}). By (\ref{14}) we get
\begin{eqnarray}
  \label{i}
\phi^{(n)}_{ \mu, \Lambda}(\zeta) & = &
\sum_{m=0}^{\infty}\frac{\zeta^m}{m!} \int_\Lambda \cdots
\int_\Lambda \varrho^{(n+m)}_\mu (x_1, \dots , x_{n+m}) dx_1 \cdots
dx_{n+m}\\[.2cm] \nonumber & = & \sum_{m=n}^{\infty} \frac{m!}{(m-n)!} (1+\zeta)^{m-n} p_\mu
(m, \Lambda) .
\end{eqnarray}
By the latter equality in (\ref{i}) one may see that $\phi^{(n)}_{
\mu, \Lambda}$ is an increasing function of $\zeta$ on $[-1,0]$,
which by the first equality yields
\begin{equation*}
\phi^{(n)}_{ \mu, \Lambda}(-1)\leq \phi^{(n)}_{ \mu, \Lambda}(0 ) =
\int_\Lambda \cdots \int_\Lambda \varrho^{(n)}_\mu (x_1, \dots ,
x_{n}) dx_1 \cdots dx_{n} \leq \left( \varkappa |\Lambda|\right)^n.
\end{equation*}
Then by (\ref{18b}) it follows that
\begin{equation}
  \label{i2}
  p_\mu (n, \Lambda) \leq \left( \varkappa |\Lambda|\right)^n/n!,
\end{equation}
that holds for any sub-Poissonian state $\mu$, where $\varkappa>0$
is its type. Note that (\ref{i2}) differs from (\ref{17}) by a
constant factor only, and hence the large $n$ asymptotic of the
occupation probabilities for sub-Poissonian states is subordinated
by the corresponding Poissonian asymptotic. If the occupation
probabilities fail to satisfy (\ref{12}), one may say that they have
a \emph{heavy tail}, cf. \cite{Kleb}. Here it is worthwhile to
recall that the states of thermal equilibrium of a particle system
with super-stable interactions are sub-Poissinian, see
\cite{Ruelle}.

\subsection{Heavy tails}

\label{SS23}

As mentioned above, if for a given state $\mu$ the large $n$ decay
of $p_\mu(n,\Lambda)$ is essentially slower than that for a Poisson
state, one says that the probabilities $\{p_\mu(n,\Lambda)\}_{n\in
\mathbf{N}}$ have a heavy tail. In particular, this means that, for
big $n$, the expected values of $n$-particle clusters
$c_\mu(n,\Lambda)$, see (\ref{6}), are much bigger than those in
sub-Poissonian states, see (\ref{9}). This can be interpreted as a
kind of \emph{clustering} occurring in such states. The aim of this
part is to give an example of the state where this effect takes
place.

The integral $\int \theta (x) \varrho(x) dx$ in (\ref{16}) can be
interpreted as the average of $\theta$ taken with respect to the
measure $\varrho(x) d x$. Let us now consider a generalization of
this average, in which it is substituted by
\[
\langle \theta \rangle = \sum_{y\in \mathbf{y}} \theta(y),
\]
where $\mathbf{y}$ is a certain infinite collection of points $y_i$,
 fixed by far. Note that this is equivalent to taking $\varrho$ in
the form $\varrho(x) = \sum_{y\in \mathbf{y}}\delta (x-y)$, where
$\delta$ is Dirac's delta-function. That is, the corresponding state
by no means is sub-Poissonian as its first correlation function is a
distribution - not a usual real valued function, which definitely
does not obey the Ruelle bound (\ref{8}). Let $\Pi_{\mathbf{y}}$ be
the corresponding state, i.e., the one for which
\[
\Pi_{\mathbf{y}}(F^\theta) = \exp\left( \sum_{y\in \mathbf{y}}
\theta (y)\right) = \prod_{y\in \mathbf{y}} e^{\theta(y)} =
F^\upsilon (\mathbf{y}), \qquad 1+\upsilon (y) = e^{\theta(y)}.
\]
Assume now that $\mathbf{y}$ is distributed over $\mathbf{R}^d$
according to the homogeneous Poisson law $\pi_\varkappa$. Then the
compound law -- Cox' cluster measure $\eta$ -- can be identified by
its values $\eta (F^\theta)$ given in the form
\begin{gather*}
  \eta(F^\theta) = \int \Pi_{\mathbf{y}}(F^\theta) \pi_\varkappa (d
  \mathbf{y}) = \int F^\upsilon (\mathbf{y}) \pi_\varkappa (d \mathbf{y}) =
  \pi_\varkappa(F^\upsilon )\\[.2cm] \nonumber = \exp\left(\varkappa\int \upsilon (x) d x
  \right) = \exp\left(\varkappa\int [e^{\theta (x)}-1] d x
  \right).
\end{gather*}
Now we take $\theta (x) = \zeta \mathbf{1}_\Lambda (x)$, see
(\ref{18a}) and (\ref{18}), and obtain
\begin{equation*}
  \phi_{\eta, \Lambda} (\zeta) = \exp\left(\varkappa\int [e^{\zeta \mathbf{1}_\Lambda (x)}-1] d x
  \right) = \exp\left(\varkappa |\Lambda|[e^{\zeta }-1]
  \right),
\end{equation*}
where $|\Lambda|$  is the volume of the vessel $\Lambda$. According
to (\ref{18b}), to get $p_\eta (n,\Lambda)$ we have to calculate the
corresponding derivatives of $\phi_{\eta, \Lambda}$. To this end, we
find the first derivative
\begin{equation*}
\phi'_{\eta, \Lambda} (\zeta) = \varkappa |\Lambda| e^\zeta
\phi_{\eta, \Lambda} (\zeta),
\end{equation*}
by which we then obtain
\begin{equation*}
\phi^{(n)}_{\eta, \Lambda} (-1) = \varkappa |\Lambda| e^{-1}
\sum_{l=0}^{n-1} \phi^{(l)}_{\eta, \Lambda} (-1).
\end{equation*}
The latter recursion can be solved in the form
\begin{gather*}
\phi^{(n)}_{\eta, \Lambda} (-1) =  a_\Lambda T_n(b_\Lambda),
\\[.2cm] \nonumber a_\Lambda = \exp\left( - \varkappa |\Lambda|(1-e^{-1})
\right), \qquad b_\Lambda = \varkappa |\Lambda| e^{-1},
\end{gather*}
where $T_n$ is Touchard's polynomial, cf. \cite[Chapt. 2]{Rio},
attributed also to J. A. Grunert and studied by S. Ranamujan, see
the historical notes in \cite[page 6]{Boy} and \cite{Berndt}. At the
same time, the first correlation function (i.e., density)
corresponding to $\eta$ is
\[
\varrho_\eta (x) = \phi'_{\eta, \Lambda}(0)/|\Lambda| = \varkappa.
\]
That is both $\pi_\varkappa$ and $\eta$ have the same densities.
However, the asymptotic properties of the corresponding occupation
probabilities are essentially different. To see this, assume that
$b_\Lambda \geq 1$, which can be obtained by taking the vessel big
enough. Then
\[
T_n(b_\Lambda) \geq b_\Lambda T_n(1) = b_\Lambda B_n,
\]
where $B_n$, $n\in \mathbf{N}$, are Bell's numbers, see
\cite{Berndt} and \cite[Chapt. 6]{Br}. Combining this with
(\ref{23}) and (\ref{18b}) we get
\begin{equation}
  \label{24a}
  p_\eta (n, \Lambda) \geq a_\Lambda b_\Lambda B_n /n!\geq
  \frac{a_\Lambda b_\Lambda}{\sqrt{2\pi} (\ln n)^n},
\end{equation}
with the latter estimate holding for big $n$. This can be compared
with, cf. (\ref{17}),
\[
p_{\pi_\varkappa} (n, \Lambda) = (\varkappa|\Lambda|)^n
e^{-\varkappa|\Lambda|} /n! \leq \frac{a_\Lambda
b_\Lambda}{\sqrt{2\pi n}}\left(\frac{[\varkappa |\Lambda|]^2}{n}
\right)^n.
\]
Thus, the occupation probabilities (\ref{24a}) decay much slower
than those for the Poisson states with the same density.

\section{The Population of Migrants}
\label{S3}

\subsection{The model}

\label{SS31}
 Now let $a^{\pm}$ and $b^{\pm}$ be the same as in
(\ref{7}). Recall that $a^{+}$ (resp. $a^{-}$) describes the
attraction (resp. repulsion) of the population members by the
existing population. At the same time $b^{\pm}$ describe the
interaction of the population members with the environment. The
interaction corresponding to $a^{-}$ can also be interpreted as
competition. Then for $\mathbf{x}$ as above, we write
\begin{equation*}
  E^{\pm}(y, \mathbf{x}) = b^{\pm} (y) + \sum_{x\in \mathbf{x}}
  a^{\pm }(y-x), \qquad y\in \mathbf{R}^d.
\end{equation*}
Thereafter, we define the corresponding gradients
\begin{eqnarray*}
  (D^{+} F)(y, \mathbf{x}) & = & E^{+} (y, \mathbf{x})\left[F(\mathbf{x}\cup y) - F(\mathbf{x})
  \right], \qquad y\in \mathbf{R}^d,\\[.2cm] \nonumber
    (D^{-} F)(y, \mathbf{x}) & = & E^{-} (y, \mathbf{x}\setminus y)\left[F(\mathbf{x}) - F(\mathbf{x}\setminus y)
  \right], \qquad y\in \mathbf{x},
\end{eqnarray*}
where $F$ is a suitable test function. Noteworthy, $(D^{+} F)(y,
\mathbf{x})$ yields the change of $F$ by adding to the population
$\mathbf{x}$ a new member located at $y$, whereas $(D^{-} F)(y,
\mathbf{x})$ corresponds to the change of $F$ by removing $y$ from
the population. Now the analog of (\ref{7d}) is the following
evolution equation, cf. (\ref{0}),
\begin{equation}
  \label{12}
  \frac{\partial}{\partial t} F_t(\mathbf{x}) = (L F_t) (\mathbf{x}) = \int (D^{+} F_t)(y,
  \mathbf{x}) dy -
  \sum_{y\in \mathbf{x}} (D^{-} F_t)(y, \mathbf{x}).
\end{equation}
In the sequel, we assume that the functions $a^{\pm}$ and  $b^{\pm}$
are continuous and bounded, i.e.
\[
0 \leq a^{\pm}(x) \leq a_{*}^{\pm}, \qquad 0 \leq b^{\pm}(x) \leq
b_{*}^{\pm}
\]
Additionally, we assume that both $a^{\pm}$ are symmetric and
integrable. That is, $a^{\pm}(x-y) = a^{\pm}(y-x)$ and
\begin{equation*}
  \int a^{\pm}(x) d x  < \infty.
\end{equation*}
The relationship between the interaction intensities $a^{\pm}$ can
be of the following types:
\begin{itemize}
 \item[(i)] long competition: there exists $\omega >0$ such that $a^{-}(x) \geq \omega a^{+}(x)$ for all $x$;
\item[(ii)] short competition: for each $\omega >0$, there exists
$x$ such that $a^{-}(x) < \omega a^{+}(x)$.
\end{itemize}
In the case of long competition, the effect of repulsion from the
population by the existing members prevails. If $b^{+}(x)\equiv 0$,
new population members appear only due to the existing ones, that
can also be interpreted as their birth. In this case, $a^{+}$ is
called \emph{dispersal kernel}. This particular case with nonzero
$a^{-}$ and $b^{-}$ is known as the Bolker-Pacala model \cite{BP1}.
In the case of long competition, $a^{-}$ decays faster than $a^{+}$,
and hence some of the attracted (by $a^{+}$) newcomers can be out of
the competitive action of the existing population. Particular cases
of the model specified by (\ref{13}) were studied in the following
works: (a) \cite{Dima,KoKoz}, case of $b^{+}(x) \equiv 0$; (b)
\cite{Koz}, case of $a^{+}(x) \equiv 0$. A detailed analysis of both
these cases in the Bolker-Pacala model based on numerical solutions
can be found in \cite{OK}.

\subsection{The statement}

\label{SS32}

 For the model considered here, the direct solution of
the Kolmogorov equation (\ref{12}) is rather impossible. Instead,
one considers its weaker version that involves states $\mu$, i.e.,
that describes their evolution $\mu_0\to \mu_t$. This is the
Fokker-Planck equation
\begin{equation}
  \label{FPE}
  \mu_t (F) = \mu_0(F) + \int_0^t \mu_s (LF) ds,
\end{equation}
which has to be solved for a sufficiently big class of test
functions $F$, for which the action of $L$ is defined. As such
functions, one can take $F^\theta$ defined in (\ref{13}) and to
restrict (\ref{FPE}) to measures satisfying (\ref{14}). To optimize
the notations we introduce variables $\mathbf{u}$ and functions
$\varrho_\mu(\mathbf{u})$ as follows. As $\mathbf{u}$ we take the
pair $\{n, (u_1, \dots , u_n)\}$ where $n\in \mathbf{N}$ and $u_i
\in \mathbf{R}^d$. Then
\begin{equation}
  \label{17}
  \varrho_\mu(\mathbf{u}) = \varrho^{(n)}_\mu (u_1, \dots , u_n) \ \
  {\rm if} \ \  \mathbf{u}= \{n,(u_1 , \dots , u_n)\},
\end{equation}
where $\varrho^{(n)}_\mu$ is the $n$-th order correlation function
of $\mu$, see (\ref{14}), which is symmetric, see (\ref{sigma}). We
extend (\ref{17}) also to the case $n=0$, i.e., to $\mathbf{u}=
\mathbf{0}:= \{0, \varnothing\}$, by setting $\varrho_\mu
(\mathbf{u})=1$. In the sequel, we call $\varrho_\mu$ the
correlation function of $\mu$ and still use `correlation function of
order $n$' if we mean $\varrho_\mu^{(n)}$. In a similar way, for a
collection of functions $g^{(n)}$, $n\geq 0$, which are symmetric as
in (\ref{sigma}), we define $g(\mathbf{u})$ by (\ref{17}). In
particular, we will deal with
\begin{equation*}
  e(\theta; \mathbf{u}) = \theta (u_1) \cdots \theta (u_n), \qquad
  e(\theta; \mathbf{0}) = 1,
\end{equation*}
with $\theta$ as in (\ref{14}). Similarly as in the case of infinite
collections $\mathbf{x}$, by writing $x\in \mathbf{u}$ we mean that
$x=u_i$ for some $i$. Likewise, for $\mathbf{u} = \{ n , (u_1 ,
\dots , u_n)\}$, we write $\mathbf{u}\cup x= \{n+1, (u_1, \dots ,
u_n, x)\}$. For $x\in
 \mathbf{u}$, i.e., for $\mathbf{u} = \{ n , (u_1 , \dots , u_{n-1}, x)\}$,
 by writing $\mathbf{u}\setminus x$ we mean $\{n-1, (u_1 , \dots ,
 u_{n-1})\}$.
With these notations, (\ref{14}) can be rewritten in the form
\begin{equation*}
  \mu (F^\theta) = \int \varrho_\mu(\mathbf{u}) e(\theta; \mathbf{u}) d
  \mathbf{u}.
\end{equation*}
The main advantage of passing to correlation functions consists in
the possibility to present $\mu(LF^\theta)$ in the form
\begin{equation*}
\mu (L F^\theta) = \int(L_\varrho \varrho_\mu)(\mathbf{u}) e(\theta;
\mathbf{u}) d  \mathbf{u},
\end{equation*}
where the operator $L_\varrho$ -- acting on correlation functions --
can be calculated from the Kolmogorov operator $L$ with the help of
a certain technique. For $L$ as in (\ref{12}), it yields
\begin{eqnarray}
  \label{21}
(L_\varrho \varrho)  (\mathbf{u}) & = & \sum_{x\in \mathbf{u}} E^{+}
(x, \mathbf{u}\setminus x) \varrho  (\mathbf{u}\setminus x)  + \int
\sum_{x\in \mathbf{u}}
a^{+} (x-y) \varrho  (\mathbf{u}\setminus x \cup y) d y \\[.2cm]
\nonumber & - & \left( \sum_{x\in \mathbf{u}}
E^{-}(x,\mathbf{u}\setminus x)\right) \varrho (\mathbf{u}) - \int
\left(\sum_{x\in \mathbf{u}} a^{-} (x-y)
\right)\varrho(\mathbf{u}\cup y) d y.
\end{eqnarray}
Employing this operator, we may pass to the following initial value
problem
\begin{equation}
  \label{22}
  \frac{\partial}{\partial t} \varrho_t = L_\varrho \varrho_t , \qquad
  \varrho_t|_{t=0} = \varrho_{\mu_0}, \ \ t>0,
\end{equation}
where $\varrho_{\mu_0}$ is the correlation function of $\mu_0$ that
appears in (\ref{FPE}). By means of the methods developed in
\cite{KoKoz,Koz} we prove the following statement.
\begin{proposition}
  \label{1pn}
Let $\mu_0$ be a sub-Poissonian state and $\varrho_{\mu_0}$ its
correlation function. Then the initial value problem in (\ref{22})
has a unique (global in time) solution $t \mapsto \varrho_t$ that
satisfies the Ruelle bound (\ref{8}) with a $t$-dependent right-hand
side and is such that, for each $t>0$, $\varrho_t$ is the
correlation function of a unique sub-Poissonian state $\mu_t$.
Moreover, these states $\mu_t$, $t\geq 0$ are such that the map
$t\mapsto \mu_t$ solves the Fokker-Planck equation (\ref{FPE}), and
hence describes the evolution of our model, in the course of which
no clustering occurs since the corresponding occupation
probabilities obey (\ref{i2}).
\end{proposition}
The mathematical proof of this statement can be done by repeating
the arguments used in the proof of the corresponding statement in
\cite{Koz}. It is  based on a nontrivial combination of various
functional-analytic methods and thus is definitely beyond the scope
of the present work.

\subsection{The evolution of correlation functions}

\label{SS33}

At first glance, besides the conclusion concerning clustering the
result formulated in Proposition \ref{1pn} has only existential
character and provides not so much information of the evolution
details. However, after further developing it may yield much more.
To see this, we begin by rewriting (\ref{22}), which -- analogously
as (\ref{7d}) -- is a symbolic (and hence concise) version of a
chain of equations, cf. (\ref{7dd}), which we are going to get now.
According to (\ref{17}), for $n\geq 2$, we have
\begin{eqnarray}
  \label{23}
(L_\varrho \varrho^{(n)}) (x_1 , \dots, x_n)& = & \sum_{i=1}^n
\left(b^{+} (x_i) + \sum_{j\neq i} a^{+} (x_i- x_j) \right)
\varrho^{(n-1)}(x_1, \dots,
x_{i-1}, x_{i+1}, \dots x_n)\qquad  \nonumber \\[.2cm] & + & \sum_{i=1}^n \int
a ^{+} (x_i - y) \varrho^{(n)} (x_1, \dots, x_{i-1}, y, x_{i+1},
\dots x_n) d y \nonumber \\[.2cm] & - & \left(\sum_{i=1}^n \left[b^{-} (x_i) + \sum_{j\neq i} a^{-} (x_i - x_j) \right] \right)
\varrho^{(n)} (x_1 , \dots, x_n) \nonumber \\[.2cm] & - &  \int
\left(\sum_{i=1}^n a^{-} (x_i-y) \right) \varrho^{(n+1)} (x_1, \dots
, x_n , y) dy.
\end{eqnarray}
The expression for $L_\varrho \varrho^{(1)}$ can be obtained from
that above by setting $\sum_{j\neq i} a^{\pm}(x_i-x_j) =0$, which is
in accord with the usual convention that $\sum_{x\in \mathbf{0}}
=0$. The latter also yields that $L_\varrho \varrho^{(0)}=0$, which
means that $\varrho^{(0)}_t = \varrho^{(0)}_0 =1$ for all $t>0$. Now
we write
\[
\frac{\partial}{\partial t} \varrho^{(n)}_t (x_1, \dots , x_n)) =
(L_\varrho \varrho_t)^{(n)} (x_1, \dots , x_n),
\]
and apply (\ref{23}) in the right-hand side of the latter. For
$n=1$, this yields
\begin{eqnarray}
  \label{24}
\frac{\partial}{\partial t} \varrho^{(1)}_t (x)& = & b^{+} (x) +
\int a^{+}(x-y) \varrho^{(1)}(y) d y \\[.2cm] \nonumber & - &
b^{-1}(x) \varrho^{(1)}(x) - \int a^{-}(x-y) \varrho^{(2)}(x,y) d y,
\end{eqnarray}
which coincides with (\ref{7}). However, now we can get the other
members of the chain.  By (\ref{23}) with $n=2$ it follows that
\begin{eqnarray}
  \label{25}
\frac{\partial}{\partial t} \varrho^{(2)}_t (x_1,x_2)& = & b^{+}
(x_1) \varrho^{(1)}(x_2) + b^{+} (x_2) \varrho_t^{(1)}(x_1) \\[.2cm] \nonumber &
+ &
\int\left[ a^{+}(x_1-y) \varrho_t^{(2)}(y,x_2) + a^{+}(x_2-y)\varrho_t^{(2)}(x_1,y)\right] d y   \\[.2cm] \nonumber & - &
\left[b^{-}(x_1) + b^{-}(x_2) + 2 a^{-} (x_1-x_2) \right]
\varrho^{(2)}(x_1, x_2)  \\[.2cm] \nonumber & - & \int \left[ a^{-}(x_1-y) + a^{-} (x_2-y) \right]
\varrho^{(3)}(x_1,x_2,y) d y.
\end{eqnarray}
Similarly, by means of (\ref{23}) one can write the next members of
the chain. By Proposition \ref{1pn} we know that the evolution of
the considered population preserves sub-Poissonicity of the system
states, which are `similar' to the Poissonian states characterized
by their densities only. Thus, the ansatz $\varrho_t^{(2)} (x_1 ,
x_2) \simeq \varrho_t^{(1)}(x_1) \varrho_t^{(1)}(x_2)$ seems
reasonable. It allows one to split the chain and to turn (\ref{24})
into (\ref{5}), which corresponds to the aforementioned mean-field
approximation, widely known in statistical physics. However, now one
can go beyond this approximation by making the ansatz in (\ref{25})
that allows for expressing $\varrho_t^{(3)}$ through
$\varrho_t^{(2)}$ and $\varrho_t^{(1)}$. The most popular one is the
Kirkwood superposition approximation \cite{BNaim}, applied in
similar situations in \cite{OK,Ruth}. If necessary, one can try to
split (\ref{22}) at higher levels.

\subsection{The kinetic equation}
\label{SS34}

In this subsection, we return to the kinetic equation (\ref{5}) and
demonstrate what kind of information can be obtained at this level
of description. Of course, first that we would have to do is proving
the existence and uniqueness of its solutions, that can be performed
similarly as in \cite[Theorem 4.2]{Dima}. Assuming this is done, we
proceed studying such solutions for various types of the model
parameters. At the beginning, we assume that the described system is
spatially homogeneous, that might be its simplest possible version.
In this case, $b^{\pm}$ -- as well as the initial density
$\varrho_0$ -- are constant functions of $x$. Correspondingly, we
will look for homogeneous solutions $\varrho_t$. This assumption
used in (\ref{5}) leads to the following nonlinear ordinary
differential equation
\begin{equation}
\label{F1} \dot{\varrho}_t = b + \alpha \varrho_t - a \varrho_t^2,
\qquad \varrho_t|_{t=0} = \varrho_0\geq 0.
\end{equation}
Here $\dot{\varrho}_t$ stands for the time derivative of $\varrho_t$
and
\begin{equation*}
  b= b^{+} \geq 0, \quad a = \int a^{-}(x) d x \geq 0, \quad \alpha = \int a^{+}
  (x) d x - b^{-}.
\end{equation*}
Note that the parameter $\alpha$ can be positive or negative,
correspondingly to the outcome of the tradeoff between the
immigration due to the attraction by the existing population and the
intrinsic emigration. Fortunately, (\ref{F1}) is a standard
differential equation of Riccati type, that can be solved
explicitly. Its solution is then presented in the form
\begin{equation}
  \label{F3}
 \varrho_t = \lambda_{+} \frac{\varrho_0 - \lambda_{-} - \delta (\lambda_{+}-\varrho_0) e^{-\omega t}}
 {\varrho_0-\lambda_{-} +(\lambda_{+}-\varrho_0)e^{-\omega t}},
\end{equation}
where
\begin{gather*}
  \lambda_{\pm} = \frac{1}{2a} \left( \alpha \pm \sqrt{\alpha^2+4ab}
  \right),  \quad  \quad \omega = \sqrt{\alpha^2 +
  4ab} , \\[.2cm] \nonumber \delta = \frac{4ab}{(\alpha+\sqrt{\alpha^2 +
  4ab})^2}.
\end{gather*}

\begin{figure}[h]
\centerline{\includegraphics[width=0.5\textwidth]{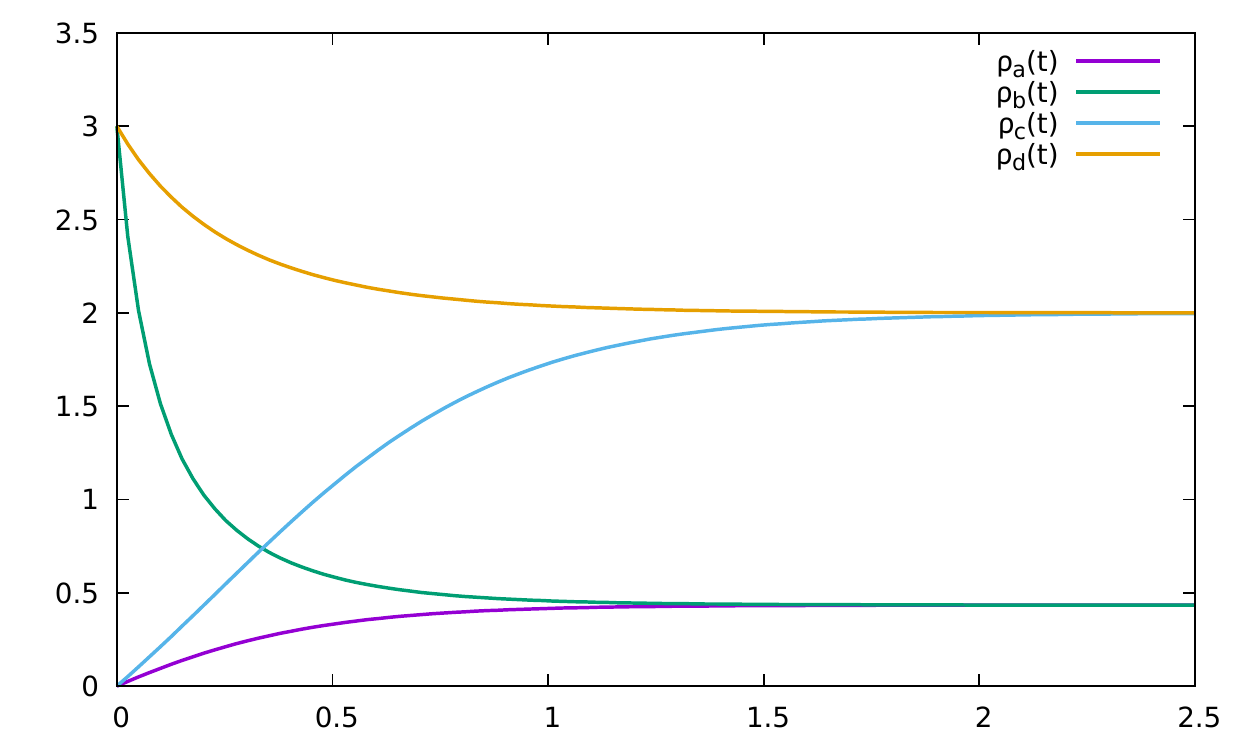}}
\caption{Solutions \eqref{F3} of equation \eqref{F1} with different
parameters and initial conditions. Value of $\varrho_t$ against $t$.
Solutions $\rho_a$ and $\rho_b$ illustrate the case with $a = 3$, $b
= 1$ and  $\alpha = -1$ with initial conditions $\varrho_0 = 0$ in
the case of $\rho_a$ and $\varrho_0 = 3$ in the case of $\rho_b$.
Solutions $\rho_c$ and $\rho_d$ correspond to the case with $a = 1$,
$b = 2$ and  $\alpha = 1$. Initial conditions as in the first case.
}\label{img_homSol}
\end{figure}

The typical behavior of the solution \eqref{F3} is plotted in Figure
\ref{img_homSol}. It has one steady state $\varrho_t=\lambda_{+}$.
Then the solution either decays or grows to the level $\lambda_{+}$,
depending on whether the initial condition exceeds it or not.

The study of the spatially homogeneous case can give a rather
superficial insight into the system evolution. Clearly, considering
non-constant parameters $b^\pm$ and/or initial condition $\rho_0$
can yield a far richer spectrum of possibilities. To illustrate some
of them, we consider the initial condition $\rho_0$ and the rates
$b^\pm$ either as constants or periodic Gaussians $G_p(c,r,s)$, see
\eqref{GaussP}. The kernels $a^\pm$ are taken as symmetrically
shifted Gaussians $G_s(c,r,s)$, see \eqref{GaussS}.

The periodic Gaussian case is defined by the formula
\begin{equation} \label{GaussP}
    G_p(x; c,r,s) = \sum\limits_{n \in \mathds{Z}} G(x - s + np; c,r),
\end{equation}
where $p$ is the period, and symmetrically shifted Gaussian by
\begin{equation} \label{GaussS}
    G_s(x; c,r,s) = \frac{1}{2} \left( G(x + s; c,r) + G(x-s; c,r) \right).
\end{equation}
Gaussian $G(c,r)$ we define as
\begin{equation*} 
    G(x; c,r) = \frac{c}{r\sqrt{2\pi}} \exp\left( - \frac{x^2}{2 r^2} \right)
\end{equation*}

In the numerical calculations presented below, the Gaussian tails
were cut after reaching small enough values and the infiniteness of
the domain was simulated using the periodic boundary conditions, see
e.g., \cite{Omel1} for a detailed description of the algorithm
analogous to that used here.

\begin{figure}[h]
\centerline{\includegraphics[width=0.5\textwidth]{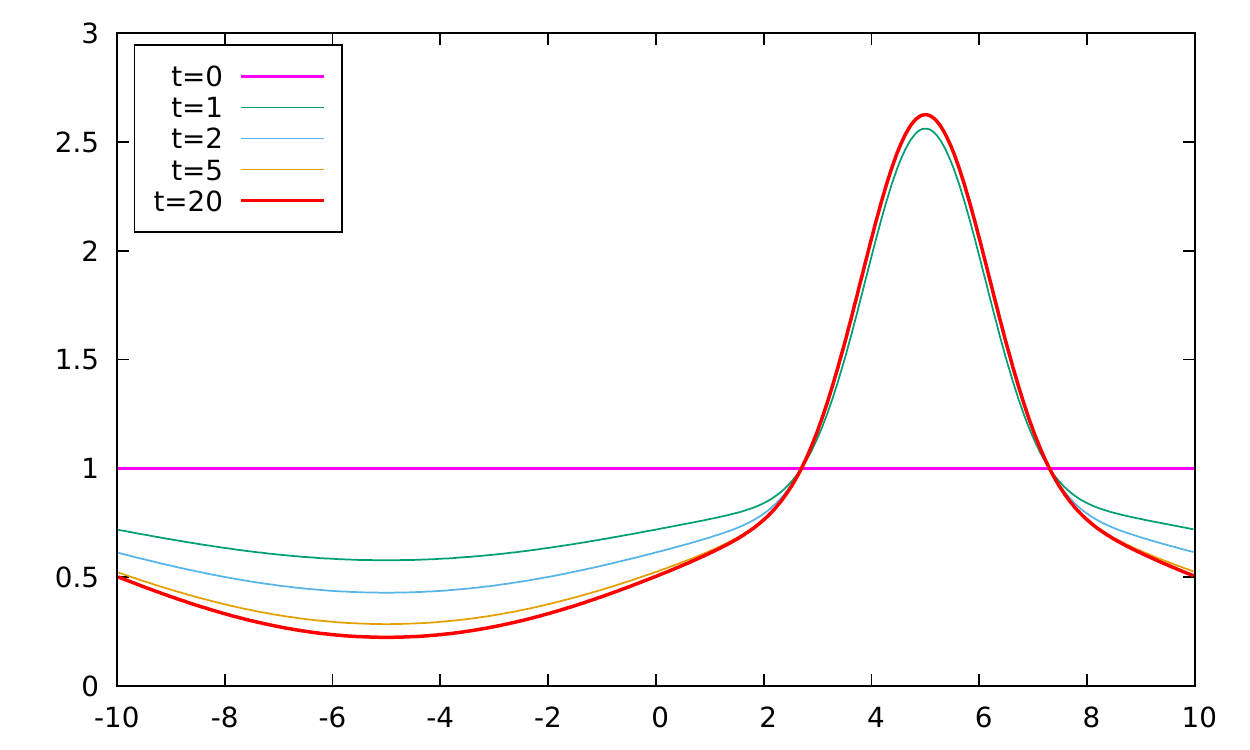} \includegraphics[width=0.5\textwidth]{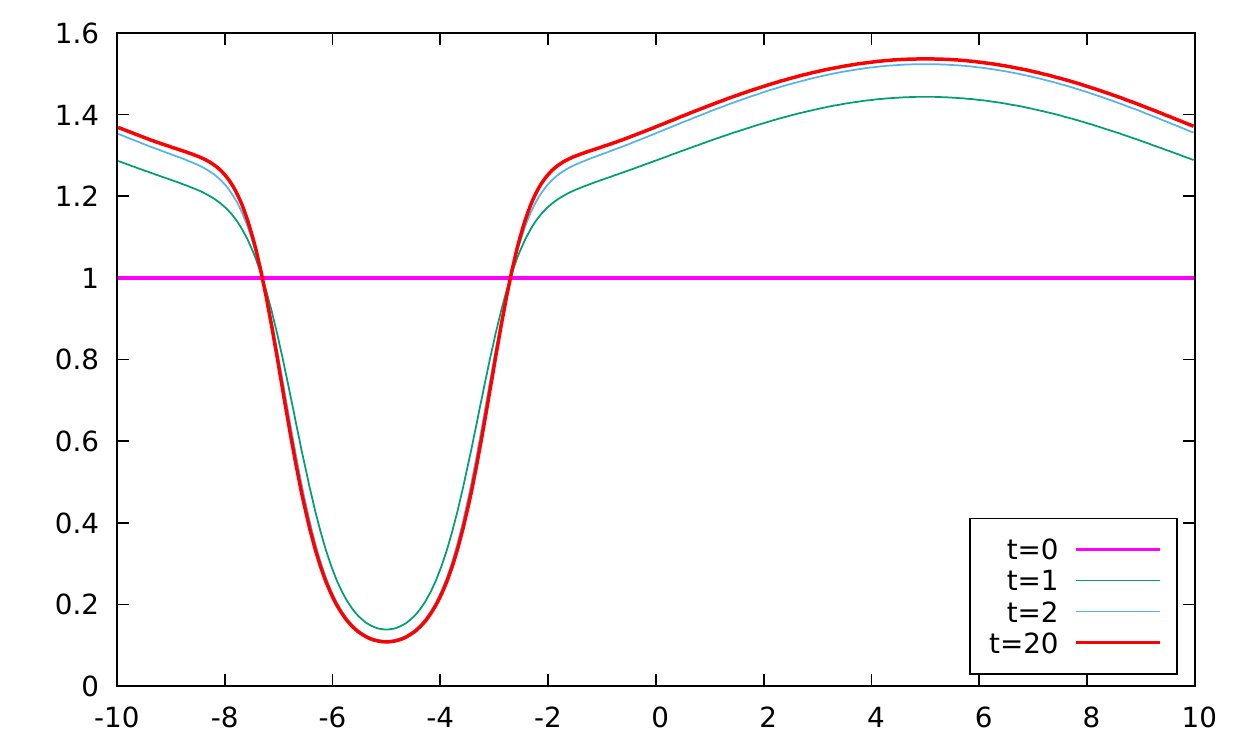}}
\caption{The role of parameters $b^+$ and $b^-$ in kinetic equation \eqref{5}. Value of $\varrho_t(x)$ versus $x$ for selected moments of time $t$. Model parameters in the form of Gaussians: $a^\pm = G(1,1)$ on both graphs, $b^+ = G_p(10,1;5)$, $b^- = G_p(10,5;-5)$ on the left graph, $b^+ = G_p(10,5;5)$, $b^- = G_p(10,1;-5)$ on the right. Constant initial condition $\rho_0 \equiv 1$ in both cases. }\label{img_parsB}
\end{figure}

In Figure \ref{img_parsB}, the role of the rates $b^\pm$ is
illustrated. Even starting from a constant initial condition, for
non-constant $b^\pm$, the system evolve towards a heterogeneous
stationary state. In the domain where $b^+$ is relatively stronger,
the density $\varrho_t$ attains higher values; it attains lower
values in the domain where $b^-$ is stronger. On the presented
graphs, one can observe how the initially constant densities
($\varrho_0 = 1$) evolve towards irregular stationary states. The
only difference between the presented cases is the range of the
rates $b^\pm$. On the left side, the range of $b^-$ is higher than
that of $b^+$,  and vice-versa on the right side.

\begin{figure}[h]
\centerline{\includegraphics[width=0.5\textwidth]{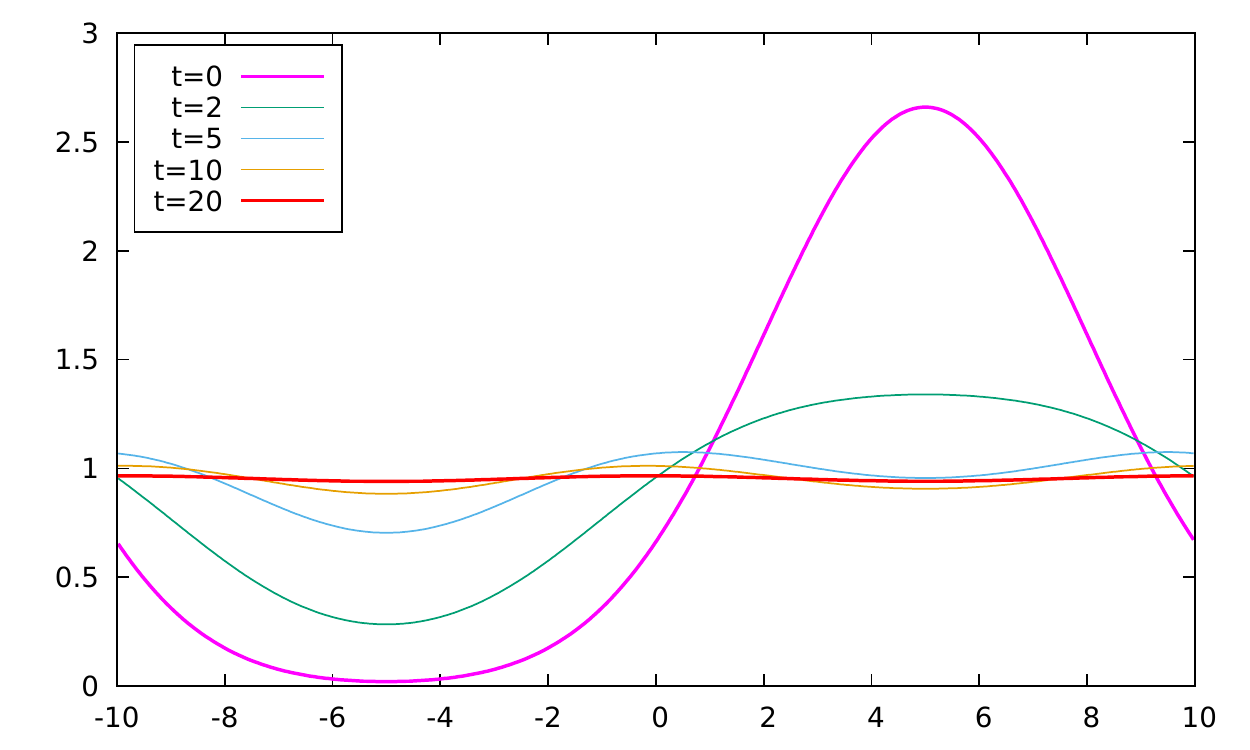} \includegraphics[width=0.5\textwidth]{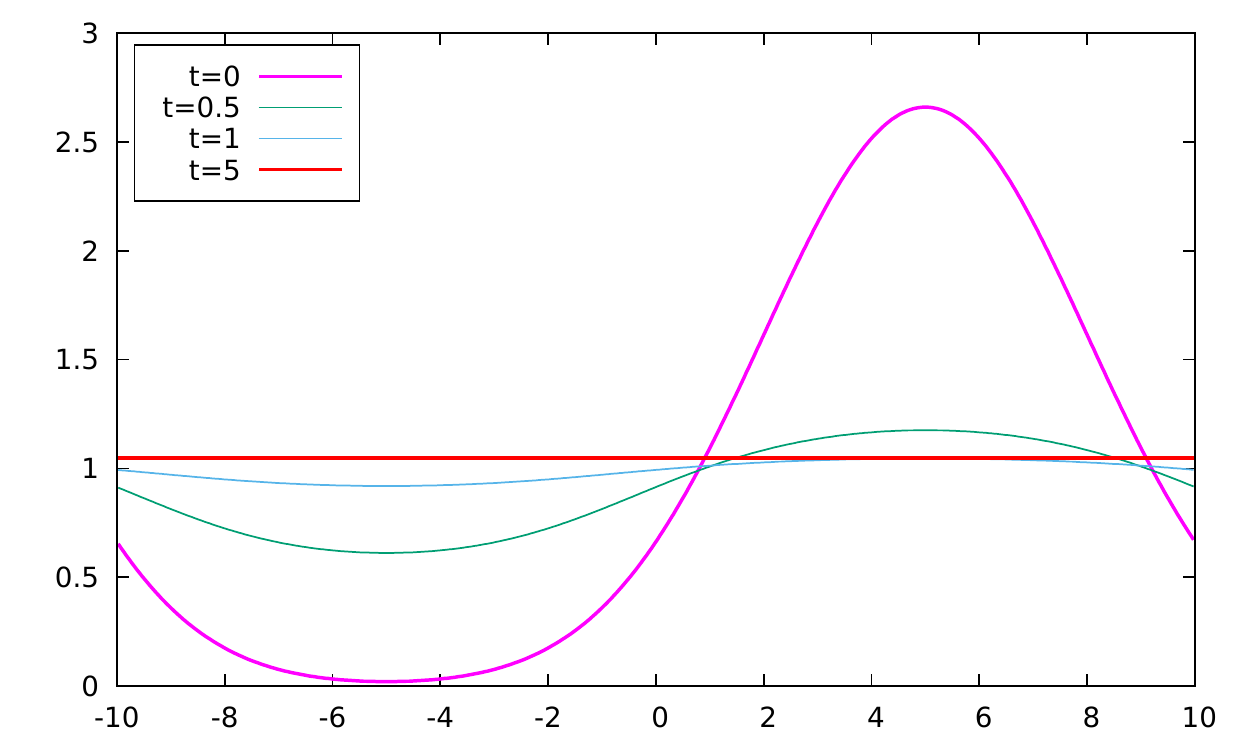}}
\caption{Short and long range of parameters $a^+$ and $a^-$ in kinetic equation \eqref{5}. Value of $\varrho_t(x)$ versus $x$ for selected moments of time $t$. Model parameters taken as $a^+ = G(2, 0.5)$, $a^- = G(2,5)$ on the left graph and $a^+ = G(2, 5)$, $a^- = G(2,0.5)$ on the right, $b^\pm \equiv 0.1$ on both. Same initial condition $\rho_0 = G_p(20,3;5)$. }\label{img_parsA}
\end{figure}

In Figure \ref{img_parsA}, the interplay of the kernels $a^\pm$ is
illustrated. In both cases, the initial condition is the same --
periodic Gaussian, and the rates $b^\pm$ are constant. The magnitude
of the kernels $a^\pm$ are comparable but they differ in the range.
It seems that the long-range $a^-$ and short-range $a^+$ (left
graph) may result in much slower approach of the stationary state
than in the opposite case (right graph).

\begin{figure}[h]
\centerline{\includegraphics[width=0.5\textwidth]{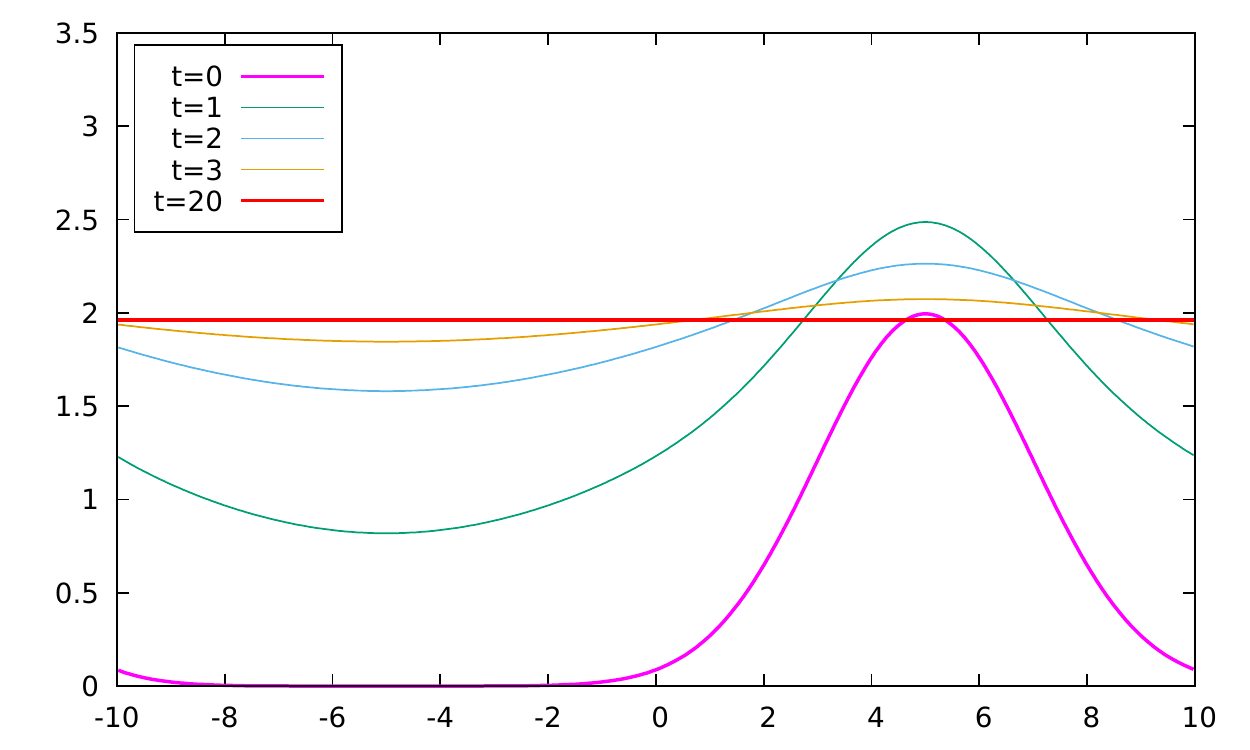} \includegraphics[width=0.5\textwidth]{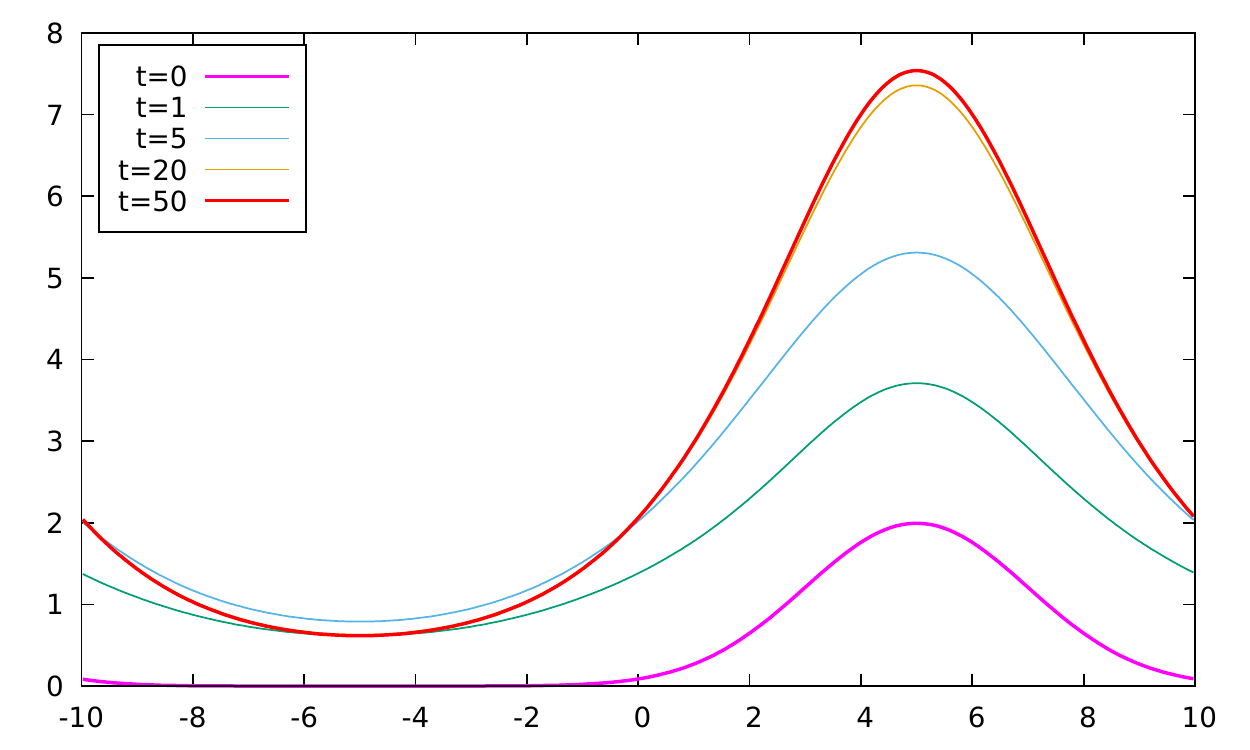}}
\caption{Shifted parameter $a^-$ in kinetic equation \eqref{5}. Value of $\varrho_t(x)$ versus $x$ for selected moments of time $t$. Model parameters $a^+ = G(2, 4)$, $b^\pm \equiv 0.1$ on both graphs, $a^-$ in the form of $G_s(1,3,5)$ on the left and $G_s(1,3;10)$ on the right. Same initial condition $\rho_0 = G_p(10,2;5)$. }\label{img_shiftedComp}
\end{figure}

\begin{figure}[h]
\centerline{\includegraphics[width=0.5\textwidth]{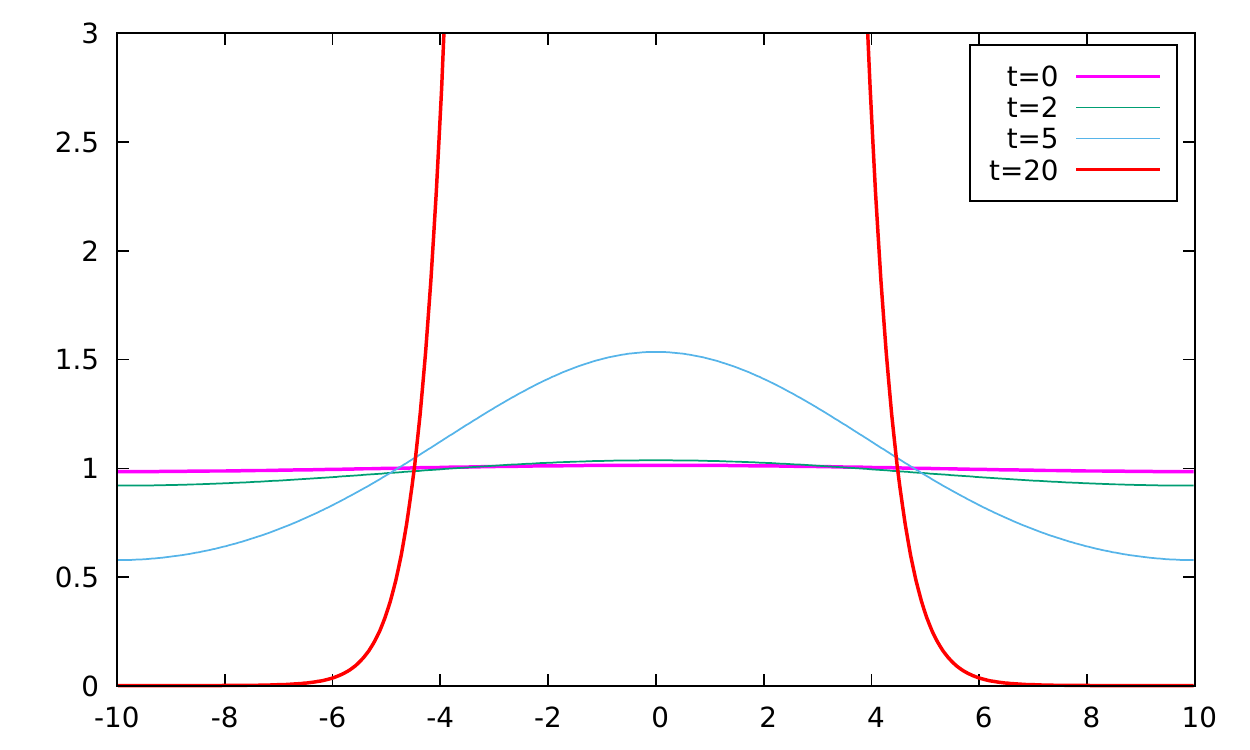}}
\caption{Unbounded growth with shifted parameter $a^-$. Value of $\varrho_t(x)$ versus $x$ for selected moments of time $t$. Model parameters $a^+ = G(1,1)$, $a^- = G_s(1,2;10)$ and $b^\pm \equiv 0.1$. Initial condition in the form of $G_p(20,10)$. }\label{img_growth}
\end{figure}

Figures \ref{img_shiftedComp} and \ref{img_growth} illustrate how
the shift of the kernels $a^-$ can substantially change the
long-time behavior of the system. In the cases presented in Figure
\ref{img_shiftedComp}, the only difference in the model parameters
is the extent of the shift of $a^-$. In the first case, where the
shift is smaller, the system evolves towards a homogeneous density,
while in the second case the bigger shift seems to result in
evolving towards a heterogeneous stationary state. When a
significant shift of $a^-$ is accompanied by a relatively strong
$a^+$ with small range, the system may experience an unbounded
growth. This is illustrated in Figure \ref{img_growth}.

More extensive numerical study of the model could reveal some other
interesting phenomena occurring in such systems, but it is beyond
the scope of this paper.

\subsection{Summary}

Let us now emphasize the basic principles of the theory developed in
this work.
\begin{itemize}

\item Similarly as in statistical physics, individual based modeling of large
 populations ought to be based on statistical/probabilistic methods
and multi-scale analysis. In particular, following Bogoliubov's idea
population states are to be described by probability measures on
appropriate phase spaces. Their evolution is obtained by solving
model-dependent Fokker-Planck equations (FPE), see (\ref{FPE}) for
our model.

\item For some models, the solutions of the Fokker-Planck
equations may have the useful property of being sub-Poissonian and
thus characterized by correlation functions. To establish this,  one
transforms the FPE into an evolution equation for correlation
functions, see (\ref{21}), (\ref{22}), solves the latter equation
and then proves that its solution is indeed the correlation function
for a solution of the FPE. The low-order correlation functions
provide such aggregated characteristics as particle density, spatial
correlations, etc.

\item Since multi-particle effects in sub-Poissonian states are inessential,
their properties are determined mostly by the low-order correlation
functions. Then the solutions of the evolution equation for
correlation functions can be approximated by solving their
`decoupled' versions (i.e., kinetic-like equations) that neglect
multi-particle effects. Noteworthy, such procedures decouple also
the macroscopic (phenomenological) theory from its microscopic
(individual-based) background. Of course, this might be admissible
only if the proof mentioned in the preceding item has been done and
also for states where multi-particle effects are inessential. At the
same time, decoupling opens the possibility to apply here powerful
methods of numerical analysis, e.g., like those developed in
\cite{Omel1,KOP,OK,Ruth}.

\end{itemize}
Here it is worth noting that, the full realization of the steps just
outlined has been done for a number of models similar to that in
(\ref{12}), see \cite{KoKoz,Koz}. For infinite systems of
interacting physical particles, this is still beyond the existing
technical possibilities, mostly because of the presence of
differential operators, see the right-hand side of (\ref{7dd}).

\section{Comparison and Discussion}
\label{SS35}

\subsection{Similar approaches}

There exists a wide variety of works dealing with kinetic equations,
or their `improved' versions based on the correlation functions of
orders two or three that take into account spatial correlations.
However, only few of them appeal to micro-states. The most typical
and, at the same time, most comprehensive recent publication of this
sort is \cite{Cornell}, where the authors deal with individual based
models defined by the Kolmogorov operator $L$, transformed
afterwards into an operator acting on correlation functions. In our
case, these operators appear in the evolution equations (\ref{12})
and (\ref{21}), respectively. In \cite{Cornell}, it has been
elaborated a general framework for calculating the time evolution of
the quantities such as our density $\varrho_t^{(1)}(x)$ and the
`truncated' correlation function $\varrho_t^{(2)}(x,y) -
\varrho_t^{(1)}(x) \varrho_t^{(1)}(y)=: u_t(x-y)$. The key step of
this approach (see \cite[page 2]{Cornell} and subsect. 1.1.4, page
8, eq. (34) of the supplementary material) consists in presenting
them  in the form of the asymptotic expansions
\begin{gather*}
\varrho^{(1)}_t (x) = q_t(x) + \varepsilon p_t(x) + o(\varepsilon),
\qquad u_t (x) = \varepsilon g_t(x) + o(\varepsilon),
\end{gather*}
where $\varepsilon>0$ is a spatial scaling parameter such that
taking the limit $\varepsilon \to 0$ corresponds to passing to a
macro-scale. By tremendous calculations (see the supplementary part
of \cite{Cornell}) the authors derive separate model-dependent
differential evolution equations  for each of $q_t$, $p_t$, $g_t$,
see \cite[eq. (2), page 3]{Cornell}, intended to cover a huge range
of models, and provide a kind of computational machine for
constructing these equations for a concrete model and then solving
them numerically and (if possible) analytically. This approach is
evidently based on the assumption that microscopic states of systems
``of an unlimited level of complexity" \cite[Abstract]{Cornell} can
  reasonably fully be specified by their low order correlation
  functions; practically, by $\varrho^{(1)}$ and $\varrho^{(2)}$.
An example where this assumption gets  problematic is provided by
Cox' cluster state described in subsect. \ref{SS23} above, for which
$\varrho^{(1)}$ is the same as for the corresponding homogeneous
Poisson state and which, therefore, cannot be distinguished from the
latter by using low order correlation functions only. At the same
time, Cox states are typical for systems of entities randomly
distributed in a fluctuating environment, which is not so exotic in
applications. In fact, instances of states with essential clustering
 -- and hence with important high-order correlations -- are well-known, see, e.g.,
\cite{Kond}. On the other hand, the very solvability of the chain
equations as in (\ref{24}), (\ref{25}) may be a problem  -- only for
few of the models mentioned in \cite{Cornell} the existence of
solutions $\varrho_t^{(n)}$ of all orders has been established, see,
e.g, \cite{Dima,KoKoz,Koz}. Moreover, the available techniques cover
only two-particle interactions, and not for all of them the
solutions exist globally in time.  Even if the solvability is
established, the corresponding solutions $\varrho^{(n)}_t$ may have
no connection to possible states of the system, and thus do not
describe the evolution of states of the considered individual based
model. To establish this connection, one has to prove  that (for
each $n$) the solution $\varrho^{(n)}_t$ is the $n$-th order
correlation function for a unique state $\mu_t$ that solves
(\ref{FPE}), and thus can be used for calculating its such relevant
characteristics as $c_{\mu_t}(k,\Lambda)$ (\ref{6}) or
$p_{\mu_t}(n,\Lambda)$, see (\ref{18b}), (\ref{18}), (\ref{14}).

\subsection{Discussion}

 The rich experience of modeling physical substances and
processes shows how important is to find a proper balance between
the model complexity and the possibility of its studying by rigorous
mathematical tools. Clearly, complex systems cannot be modeled in an
adequate way by simple models. On the other hand, for complex models
the mathematically rigorous results of their study are usually too
general and thus tell almost nothing to the modelers. The typical
way of finding balance here is to employ usually uncontrolled
approximations, often called `ansatz'. The mean field approximation
mentioned in subsect. 2.1 above is among the most popular. Such
approximations yield a more detailed description, whereas the
possible weakness of their mathematical justification may be
compensated by comparison with the corresponding experimental data.
Fortunately for physicists, their experiments are incomparably more
reproducible than those in other natural sciences. To their
consolation, mathematicians involved in modeling physical systems
are able to show (at least in some cases) that a given ansatz yields
an exact result but for another model, that is not purely
individual-based and should be considered as a toy model -- a
convenient `caricature' of the initial model. As good examples may
serve individual based models with Curie-Weiss (nonlocal)
interactions, see, e.g., \cite{KK}, providing a `mathematical
justification' of the naive mean field approximation.

Since the experimental results concerning the majority of the
complex systems studied in life sciences are much less reliable than
those mentioned above, the accuracy of the theoretical modeling of
such systems predetermines the correctness of their understanding.
Hence, one can safely say that the adequate description of
individual based models ``containing interactions of an unlimited
level of complexity" -- either by means of the technique of
\cite{Cornell}, or by any other machinery of this kind -- is barely
possible.

The principal novel aspects of the present study can be summarized
the following.
\begin{itemize}
\item The complete utilization of an individual based model
characterized by its Kolmogrov operator $L$, should be performed by
constructing the evolution of states that solve the corresponding
Fokker-Planck equation. Such solutions provide complete information
concerning the evolution of the considered model. To facilitate this
study, one can use a suitable class of states to which the possible
solutions of the latter are confined. The most natural one is the
class of sub-Poissonian states defined in subsect. \ref{SS22}, the
members of which are characterized by their correlations functions.
In such states, clustering is not essential and thus the spatial
scaling procedures developed and used in \cite{Cornell}
\underline{minimally} distort the description. Significantly,
relevant states of systems of interacting  physical particles are
sub-Poissonian \cite{Ruelle}. At the same time, many of states
important in applications are not sub-Poissonian, and may not have
correlation functions at all.
\item One can transform the
Fokker-Planck equation to an evolution equation solutions of which
may (but need not) be the correlation functions for some
sub-Poissonian states. Such equations appear as \cite[eq. (13), page
6 od Suppl.]{Cornell} and (\ref{22}) in our work. Then one has to
prove that their solutions are indeed correlation functions, which
is missing in the framework of \cite{Cornell}. This step establishes
a crucial \underline{link} between the evolution of the initial
individual based model and its description in terms of
$\varrho_t^{(n )}$. Without it, existence of such links is nothing
more than wishful thinking.
\end{itemize}


\begin{thebibliography}{00}


\bibitem{Neuhauser}  C. Neuhauser, Mathematical challenges in spatial
ecology, Notices of AMS, 48 (11) (2001) 1304--1314.

\bibitem{BP3}  B. M. Bolker, S. W.  Pacala, C. Neuhauser,
 Spatial dynamics in model plant communities: What do we really
know?  The American Naturalist,  162 (2003) 135--148.



\bibitem{Burini} D. Burini, L. Gibelli, N. Outada,  A kinetic theory approach
to the modeling of complex living systems. In: Bellomo N., Degond
P., Tadmor E. (eds) Active Particles, Volume 1. Modeling and
Simulation in Science, Engineering and Technology. Birkh\"auser,
Cham. 2017.


\bibitem{otso} O. Ovaskainen, B. Meerson, Stochastic models of population
extinction, Trends in Ecology $\&$ Evolution 25 (2010), 643--652.

\bibitem{Cornell} S. J. Cornell, Y. F. Suprunenko, D. Finkelshtein, P. Somervuo, O. Ovaskainen,
A unified framework for analysis of individual-based models in
ecology and beyond, Nat. Commun. 10 (2019), 4716 1--14.

\bibitem{OSF} O. Ovaskainen, P. Somervuo, D. Finkelshtein, A general
mathematical method for predicting spattio-temporal correlations
emerging from agen-based models, J. R. Soc. Interface, 17 (2020),
20200655 1--10.

\bibitem{Newman} M. E. J. Newman, Resource letter CS-1: complex
systems, Am. J. Phys. 79 (2011), 800--810.


\bibitem{Green}, D. G. Green, S. Sadedin, Interactions matter -- complexity in landscapes
and ecosystems, Ecological Complexity 2 (2005) 117--130.

\bibitem{DeA} D. L. De Angelis, W. M. Mooij, Individual-based modeling of ecological and evolutionary processes,
Annu. Rev. Ecol. Evol. Syst. 36 (2005) 147--168.

\bibitem{KozF} Yu. Kozitsky,  Dynamics of spatial logistic model: finite systems, in Semigroups of Operators -Theory and Applications,
J.  Banasiak, A. Bobrowski, M. Lachowicz, eds, Springer Proceedings
in Mathematics $\&$ Statistics, vol 113. Springer, Cham, 2015.

\bibitem{Simon}  B. Simon, The Statistical Mechanics of Lattice
Gases. I, Princeton University Press, Pronceton, NJ, 1993.


\bibitem{Bogol} N. N. Bogoliubov, Problems of Dynamical Theory in Statistical Physics
[Russian], OGIZ, Gostechnizdat, Moscow-Leningrad, 1946.

\bibitem{DSS} R. L. Dobrushin, Y. G. Sinai, Y. M. Sukhov, Dynamical systems of statistical
mechanics, in Dynamical Systems II. Encyclopaedia of Mathematical
Sciences, vol 2., Y. G. Sinai, eds, Springer, Berlin, Heidelberg,
1989.




\bibitem{Dima} D. Finkelshtein, Yu. Kondratiev, Yu. Kozitsky, O.
Kutoviy, The statistical dynamics of a spatial logistic model and
the related kinetic equation, Math. Models Methods Appl. Sci. 25
(2015) 343--370.

\bibitem{KoKoz} Yu. Kondratiev, Yu. Kozitsky, The evolution of
states in a spatial population model, J. Dyn. Diff. Equat., 30
(2018) 135--173.

\bibitem{Koz}  Yu. Kozitsky,  Evolution of infinite populations of immigrants: micro- and mesoscopic
description, J. Math. Anal. Appl.,  477 (2019) 294--333.


\bibitem{OK} I. Omelyan, Yu. Kozitsky, Spatially inhomogeneous population dynamics: beyond the mean field
approximation, J. Physics A Mathematical and Theoretical 52 (2019)
305601

\bibitem{Omel1}  I. Omelyan, Yu. Kozitsky, K. Pilorz,
 Algorithm for numerical solutions to the kinetic equation of a
spatial population dynamics model with coalescence and repulsive
jumps,  Numer. Algor. (2020),
https://doi.org/10.1007/s11075-020-00992-9


\bibitem{KOP}  Yu. Kozitsky, I. Omelyan, K. Pilorz, Jumps and coalescence in the continuum: A numerical
study,  Appl. Math. Comput., 390 (2021) 125610


\bibitem{Ruelle} D. Ruelle, Superstable interactions in classical statistical mechanics,  Comm. Math.
Phys.,  18 (1970) 127--159.

\bibitem{Rue} D. Ruelle, Existence of a phase transition in a continuous classical system, Phys.
Rev. Lett.,   27 (1971) 1040--1041.


\bibitem{BP1} B. M. Bolker, S. W. Pacala,
 Using moment equations to understand stochastically driven spatial
pattern formation in ecological systems,  Theoret. Population Biol.
52 (1997) 179--197.



\bibitem{KK} Yu. Kozitsky, M. Kozlovskii, A phase transition in a
Widom-Rowlinson model with Curie-Weiss interaction, {\it J. Stat.
Mech. Theory Exp.} 2018, no. 7, 073202, 28 pp.


\bibitem{Lenard}  A. Lenard, Correletion functions and the uniqueness
of the state in classical statistical mechanics,  Comm. Math. Phys.
30 (1973) 35--44.



\bibitem{Mu}  D. J.  Murrell, U.  Dieckmann, R.  Law, On moment closures for population
dynamics in contunuous space,   J. Theoret. Biol.   229 (2004)
421--432.




\bibitem{Plank} M. J. Plank, M. J. Simpson, R. N. Binny, Small-scale spatial
structure influences large-scale invasion rates. Theor Ecol 13
(2020) 277--288.


\bibitem{ST} J. A. Shohat, J. D. Tamarkin, The Problem of Moments. Amer. Math.
Soc., Providence, R.I., 1943.



\bibitem{DV} D. J. Daley, D. Vere-Jones, An Introduction to the
Theory of Point Processes. Vol. I. Elementary Theory and Methods.
Second edition. Probability and its Applications (New York),
Springer-Verlag New York, 2003.

\bibitem{Kleb} V. Klebanov, Heavy Tailed Distributions, Matfyzpress,
Prague, 2003.

\bibitem{Rio} J. Riordan, Combinatorial Identities, John Wiley $\&$
Sons, Inc., New York -- London -- Sidney, 1968.

\bibitem{Boy} K. N. Boyadzhiev, Exponential polynomials, Stirling
numbers, and evaluation of some gamma integrals. Abstr. Appl.
Analysis Volume 2009, article ID 168672, 18 pp.

\bibitem{Berndt} B. C. Berndt, Ramanujan reaches his hand from his
grave to snatch your theorems from you, Asia Pac. Math. Newsl. 1
(2011) 8--13.

\bibitem{Br} N. G. de Bruijn, Asymptotic Methods in Analysis, Third
edition, Dover Publishers, Inc., New York, 1981.


\bibitem{Kond} D. L. Finkelshtein, Yu. G. Kondratiev,  O. Kutoviy, Individual based
model with competition in spatial ecology, SIAM J. Math. Anal. 41
(2009) 297--317.



\bibitem{BNaim} A. Ben-Naim, The Kirkwood superposition approximation, revisited and reexamined,  J. Adv.
Chem. 1 (2013) 27--35.

\bibitem{Ruth} R. E. Baker, M. J. Simpson, Correcting mean-field approximations for
birth-death-movement processes, Phys. Rev. E  82 (2010) 041905.









\end{thebibliography}
\end{document}